%% file: arxiv_v2.tex
\documentclass[hidelinks,onefignum,onetabnum]{siamart250211}

\input{ex_shared}

\ifpdf
\hypersetup{
  pdftitle={An Example Article},
  pdfauthor={D. Doe, P. T. Frank, and J. E. Smith}
}
\fi

\makeatletter
\g@addto@macro\normalsize{%
  \setlength{\abovedisplayskip}{2pt}%
  \setlength{\belowdisplayskip}{3pt}%
  \setlength{\abovedisplayshortskip}{2pt}%
  \setlength{\belowdisplayshortskip}{3pt}%
}
\makeatother

\headers{A parallel framework for graphical optimal transport}{J. Fan, I. Haasler, Q. Zhang, J. Karlsson, Y. Chen}

\title{A parallel framework for graphical optimal transport \thanks{ 
\funding{This work was supported by the Knut and Alice Wallenberg
foundation under grant KAW 2021.0274, the Swedish Research Council (VR) under grant 2020-03454, and the NSF under grant 2206576.}}}

\author{Jiaojiao Fan\thanks{NVIDIA,
\email{jiaojiaof@nvidia.com}. The first two authors contributed equally.}
       \and
 Isabel Haasler\thanks{Department of Information Technology, Uppsala University, Uppsala, Sweden, \email{isabel.haasler@it.uu.se }}
       \and
Qinsheng Zhang\thanks{NVIDIA, \email{qinshengz@nvidia.com} }
       \and
Johan Karlsson\thanks{Department of Mathematics,
       KTH Royal Institute of Technology,
       Stockholm, Sweden,
 \email{johan.karlsson@math.kth.se} }
       \and
Yongxin Chen\thanks{School of Aerospace Engineering,
       Georgia Institute of Technology,
       Atlanta, Georgia, USA   \email{yongchen@gatech.edu} }
       }

\begin{document}

\maketitle

\begin{abstract}
We study multi-marginal optimal transport (MOT) problems where the underlying cost has a graphical structure. These graphical multi-marginal optimal transport problems have found applications in several domains including traffic flow control,
barycenter
and regression problems in the Wasserstein space, and
Hidden Markov model inference problems.
The MOT problem can be approached through two formulations: a single big MOT problem, or coupled minor OT problems. In this paper, we focus on the latter approach and demonstrate its efficiency gain from parallelization. For tree-structured MOT problems, we introduce a novel parallelizable algorithm that significantly reduces computational complexity. Additionally, we adapt this algorithm for general graphs, employing the modified junction trees to enable parallel updates. Our contributions, validated through numerical experiments, offer new avenues for MOT applications and establish benchmarks in computational efficiency.
\end{abstract}



\section{Introduction}


Optimal transport (OT) theory offers a fundamental mathematical framework for the comparison of probability distributions~\cite{peyre2017computational}. Central to OT is the quest for the most cost-efficient strategy to morph one distribution into another. With significant strides in the development of computationally efficient methods~\cite{Cut13, zanetti2023interior}, OT has burgeoned into a plethora of applications. It has made substantial impacts in control theory \cite{CheGeoPav14e,CheGeoPav21a,SinHaaZha20,ringh2021efficient}, generative modeling \cite{fan2021neural,fan2021variational,korotin2022neural}, image processing \cite{peyre2017computational,fan2021neural,zhou2022efficient}, domain adaptation \cite{seguy2017large,perrot2016mapping}, Bayesian statistics \cite{fan2020scalable,srivastava2015wasp}, and biology \cite{bunne2022proximal,demetci2020gromov}, showcasing the versatility and the expansive scope of OT as a tool for analytical and computational advancement.
Multi-marginal optimal transport (MOT)~\cite{ba2022accelerating, beier2023unbalanced, trillos2023multimarginal, yang2023estimating}, a generalization of OT, involves optimizing over transports among more than two distributions, thereby extending the applicability of OT to a broader range of problems~\cite{Pas12,Pas15,Nen16}. Its applications are vast and varied, encompassing areas such as estimation and control~\cite{ringh2023mean,ringh2021efficient}, signal processing~\cite{ElvHaaJakKar20}, and machine learning \cite{HaaSinZha20, noble2023tree}.

The MOT problem frequently embodies a graphical structure where each node denotes a distribution, and the edges represent the transport plans. 
For example, the well-known Wasserstein barycenter problem finds a generalized average of a set of distributions \cite{AguCar11}. This is done by finding a distribution that minimizes the sum of optimal transport distances to the given distributions, which can be represented by a star-shaped graph \cite{fan2021complexity}.
Originally, the Wasserstein barycenter problem was formulated as a sum of coupled bi-marginal OT problems \cite{carlier2010matching, rabin2012wasserstein}. However, it can equivalently be posed as a MOT problem, where the graphical structure is instead invoked through the underlying cost.
This equivalence between the MOT formulation and the coupled OT formulation holds for any structured problems, where the underlying graph is a tree, i.e., an acyclic graph \cite{HaaRinChe20}.
Although the coupled OT methodology often presents a more intuitive formulation, its computational development has not been extensively pursued, and existing numerical solvers for general graphical structures rely predominantly on the MOT formulation \cite{ElvHaaJakKar20, HaaRinChe20, HaaSinZha20, fan2021complexity}.
To close this gap in the literature, we develop numerical methods for graph-structured OT problems that rely on the coupled OT formulation.

More precisely, following the powerful standard approach of adding an entropic term to the OT objective \cite{Cut13}, we first develop a novel scheme for a sum of OT problems that are coupled according to a tree-structure.
We also prove complexity results for solving the problem to a desired accuracy.
It turns out that this algorithm can be easily parallelized, thus significantly improving the runtime of previous algorithms that are based on the equivalent MOT problem.
We also extend our algorithm to the setting of general graph-structures that may contain cycles. In this setting, the algorithm is based on factorizing the underlying graph as a modified junction tree and solving an OT problem that is structured according to this modified junction tree.
Thus, in this case the problem is posed as a sum of 
several coupled MOT problems of small size, where each MOT problem corresponds to a clique, i.e., a collection of nodes, in the modified junction tree. To summarize, our contributions are the following.


(i) For MOT problems associated with tree structures, we reformulate them as coupled bi-marginal OT problems and introduce a novel algorithm for solving them with local entropy regularization in \cref{sec:tree_algo}. In our approach, we represent the tree as a bipartite graph and categorize the optimization variables into two groups. This partitioning facilitates the parallelization of our algorithm. This advancement bridges the existing gap in the literature by providing a dedicated algorithmic solution for coupled small OT problems.

(ii) In \cref{sec:tree_complexity}, we 
{conduct} 
a thorough complexity analysis for our algorithm addressing coupled bi-marginal OT problems within tree structures. Contrasting with standard MOT problems that iterate over nodes sequentially, our parallelizable algorithm requires only $\widetilde\cO \left( {|E|^2}/{\delta^2} \right)$ iterations. 
Here, $|E|$ denotes the total number of edges,
and $\delta$ specifies the desired accuracy. This represents a significant efficiency gain, reducing the number of iterations by a factor of $|E|$ compared to the standard MOT problem approach analyzed in \cite{fan2021complexity}.


(iii) In \cref{sec:extension}, 
we introduce the notion of a modified junction tree, allowing us to extend the scope of our algorithm from tree structures to general graphs. This extension involves formulating the problem as a sum of coupled small MOT problems and constructing a bipartite graph to enable a parallel algorithm.
To validate the efficacy and robustness of our algorithm, we have conducted several numerical experiments, the results of which affirm the algorithm's applicability to graph structures.

\section{Preliminaries and problem formulation}\label{sec:pre}

In this section, we review some background material and introduce the problem we consider in this work.
In particular, we first define optimal transport and multi-marginal optimal transport and recap the iterative scaling algorithm for approximately solving these problems.
We then introduce MOT problems with graphical structure and provide several motivating examples for considering these types of problems.

\subsection{Optimal transport and iterative scaling}
\subsubsection{Optimal transport} \vspace{-5pt}
The OT problem seeks an optimal transport plan that moves mass from a source distribution to a target one with minimum cost.
In this work, we focus on the discrete OT problem (see, e.g., \cite{Vil03} for the continuous formulations).
Thus, consider two non-negative vectors $\mu_1\in \mR_+^{d}$ and $\mu_2 \in \mR_+^{d}$  
with mass $1$, i.e., $\sum_{x_1} \mu_1(x_1) = \sum_{x_2} \mu_2 (x_2) = 1$, that describe discrete probability distributions. For simplicity of notation, we assume the marginals have equal dimension $d$ throughout, though the extension to the setting with unequal dimensions is straightforward.
Here $\mu_1(x)$ and $\mu_2(x)$ denote the amount of mass in the source distribution and target distribution, respectively, at location $x$.
We also define a cost matrix $C\in \mR^{d\times d}$, where $C(x_1,x_2)$ describes the transport cost of moving a unit mass from point $x_1$ to $x_2$.
Similarly, we describe the transport by a non-negative matrix $B \in \mR_+^{d\times d}$, where the element $B(x_1,x_2)$ describes the amount of mass transported from point $x_1$ to $x_2$.
This defines the optimal transport problem
\begin{equation}\label{eq:omt_bi_discrete}
\begin{aligned}
\min_{B \in \mR_+^{d\times d}}~ & \tr( C^\top B) 
\qquad \text{ subject to }    & B \ett = 
\mu_1, \quad B^\top \ett = 
\mu_2,
\end{aligned}
\end{equation}
where $\ett$ denotes a vector of ones of proper dimension. The constraints enforce that $B$ is a feasible transport plan between $\mu_1$ and $\mu_2$, and the objective function 
$\tr(C^\top B) = \sum_{x_1,x_2} C(x_1,x_2) B(x_1,x_2)$ describes the total cost of transportation for the given transport plan. 

\subsubsection{Multi-marginal optimal transport} \vspace{-5pt}
The MOT problem extends the OT problem \cref{eq:omt_bi_discrete} to finding a transport plan between multiple distributions $\mu_1,\dots,\mu_J$, where $J>2$.
Let $\mu_j \in \mR^{d}$ for $j=1,\dots,J$. Then the multi-marginal transport cost and transport plan are described by $J$-mode tensors $\bC \in \mR^{d\times \dots \times d}$ and $\bB \in \mR_+^{d\times \dots \times d}$. 
Here, $\bC(x_1,x_2,\ldots, x_J)$ and $\bB(x_1,x_2,\ldots, x_J)$ denote the cost and amount of transport associated with the locations $x_1, x_2,\ldots, x_J$.
Thus, analogously to \cref{eq:omt_bi_discrete}, the MOT problem is defined as
\begin{equation} \label{eq:omt_multi_discrete}
		\min_{\bB \in \mR_+^{d\times \dots \times d}}  \langle \bC, \bB \rangle                                   
	\qquad	\text{ subject to }                                P_j (\bB) = \mu_j,  \text { for } j \in \Gamma,
\end{equation}
where $\Gamma\subset \{1,2,\dots,J\}$ is an index set specifying which marginal distributions are given, and the projection on the $j$-th marginal of $\bB$ is computed as
\begin{equation} \label{eq:proj_discrete}
	P_j(\bB) = \sum_{x_1,\dots,x_{j-1},x_{j+1},x_J} \bB (x_1,\dots,x_{j-1},x_j,x_{j+1},\dots,x_J).
\end{equation}
In the original MOT formulation \cite{GanSwi98,Pas15}, constraints are given on all marginal distributions, that is the index set in the constraints is $\Gamma = \{1,2,\dots,J\}$.
However, as we will see in the next section, in many applications not all marginal distributions are given, thus in this work $\Gamma$ is typically a strict subset of $\{1,2,\dots,J\}$.
Finally, we note that the standard bi-marginal OT problem \cref{eq:omt_bi_discrete} is a special case of the MOT problem \cref{eq:omt_multi_discrete} with $J=2$ and $\Gamma=\{1,2\}$.

\subsubsection{Iterative scaling algorithm} \vspace{-5pt}
We note that MOT is a standard linear program. However, as it has 
$d^J$
variables, solving it directly scales exponentially with the number of marginals $J$.
The computational burden of solving \cref{eq:omt_multi_discrete} can be partly alleviated by regularizing the objective with an entropy term.
This approach has first been proposed for the bi-marginal problem \cref{eq:omt_bi_discrete} in  \cite{Cut13,CheGeoPav21b}, and has later been adapted to multi-marginal problems \cite{BenCarCut15, HaaRinChe20, HaaSinZha20}.
The entropy regularized MOT problem is the strictly convex optimization problem
\begin{equation} \label{eq:omt_multi_regularized}
	\begin{aligned}
		\min_{\bB \in \mR^{d\times \dots \times d }}  \langle \bC, \bB \rangle + \epsilon  \cH(\bB | \bM)  
	\qquad	\text{ subject to }                              P_j (\bB) = \mu_j,  \text { for } j \in \Gamma,
	\end{aligned}
\end{equation}
where $\epsilon>0$ is a small regularization parameter and 
%
\begin{equation} \label{eq:entropy_term}
 \cH(\bB | \bM)
 = \!\!\!  \sum_{x_1,\dots,x_J} \!\!\!  \bB(x_1,\dots,x_J)  \left( \log \left( \bB(x_1,\dots,x_J) \right) - \log \left( \bM(x_1,\dots,x_J) \right) -1 \right) \!\!
\end{equation}
with $\bM(x_1,\dots,x_J) = \prod_{j\in\Gamma} \mu_j(x_j)$.
The entropy regularized problem \cref{eq:omt_multi_regularized} can be solved by an iterative scaling algorithm, presented in \cref{alg:sinkhorn}.
This method is also known as Sinkhorn iterations, especially in the bi-marginal setting \cref{eq:omt_bi_discrete}.

\begin{algorithm*}[tb]
\caption{Iterative Scaling Algorithm for MOT 
}
\label{alg:sinkhorn}
    Compute $\bK=\exp(- \bC/\epsilon)$, 
    initialize $u_1, u_2, \ldots, u_J$ 
    to $\mathbf{1}$

    \While{not converged}{
    \For{$j \in \Gamma$} {
    Update  $\bU \leftarrow u_1 \otimes u_2 \otimes \dots \otimes u_J$
    
    Update $u_j \leftarrow u_j \odot \mu_j ./ P_j(\bK \odot \bU \odot \bM )$
    }
		}
  \Return{$\bB = \bK \odot \bU \odot \bM $}
\end{algorithm*}
%

In the following, we briefly sketch how to derive 
 \cref{alg:sinkhorn} as a block coordinate ascent method in the dual of \cref{eq:omt_multi_regularized}, for details we refer the reader to \cite{ElvHaaJakKar20,HaaRinChe20,HaaSinZha20}.
Namely, utilizing Lagrangian duality theory, one can show that the optimal solution to \cref{eq:omt_multi_regularized} is of the form $\bB = \bK \odot \bU  \odot \bM$,
where $\odot$ denotes element-wise multiplication and the tensors are given by
\begin{equation}
	\bK  = \exp(- \bC/\epsilon,) \quad
 	\bU  = u_1 \otimes u_2 \otimes \dots \otimes u_J,
    \quad 
    u_j = \begin{cases} \exp\left( \frac{\lambda_j}{\epsilon} \right), & \text{ if } j \in \Gamma \\
		 \ett,                                    & \text{ otherwise,}\end{cases}
    \label{eq:U}
\end{equation}
where $\lambda_j\in \mathbb{R}^{d}$ is the dual variable corresponding to the constraint $P_j (\bB) = \mu_j$ on the $j$-th marginal.
The corresponding dual problem of \cref{eq:omt_multi_regularized} is
\begin{equation} \label{eq:multi_omt_dual}
	\max_{\{\lambda_j, j\in \Gamma\}}  -\epsilon \langle\bK, \bU \rangle {+ \sum_{j \in \Gamma} \lambda_j^\top \mu_j,}
\end{equation}
where $\bU$ depends on the dual variables as in \cref{eq:U}.
 \Cref{alg:sinkhorn} is then derived as a block coordinates ascent in this dual and is therefore guaranteed to converge globally.

\subsection{Graph-structured optimal transport}

In this work we consider multi-marginal optimal transport problems with a cost that is composed as a sum of cost terms, which each depend on a subset of the marginals.
That is, we consider cost tensors of the form 
\begin{equation}\label{eq:cost_structure}
	\bC(\bx) = \sum_{\alpha\in F} C_{\alpha}(\bx_\alpha),
\end{equation}
where $F$ is a set of subsets of the tensors indices, i.e., $\alpha \subset  \{1,\dots,J\} $ for each element $\alpha \in F$.
If we associate the index set of the transport tensor with a set of vertices, $V = \{1,\dots,J\}$, then the set $F$ can be understood as the set of factors in a factor graph \cite{HaaSinZha20, fan2021complexity}.
In particular, all subsets $\alpha \in F$ that contain two elements can be illustrated as edges in a graph.
Moreover, we allow for constraints that act on projections onto several marginals of the transport tensor $\bB$.
For a set $\alpha \subset V$, we define the projection
 $   P_\alpha (\bB) = \sum_{x_j,\  j \in V \backslash \alpha} \bB(x_1,\dots, x_J).$
Let $F_\Gamma \subset F$ denote the set of factors which are constrained, and for $\alpha \in F_\Gamma$, let $\boldsymbol\mu_\alpha \in \mR_+^{d^{|\alpha|}}  $  be the given tensor, 
where the number of modes equals the cardinality $|\alpha|$ of $\alpha$.
This let us define the graph-structured MOT problem as 
\begin{equation} \label{eq:ot_graph}
	\begin{aligned}
		\min_{\bB \in \mR_+^{
  d \times \dots \times  d 
  }}  \langle \bC, \bB \rangle                                                   
	\qquad	\text{ subject to }                                P_\alpha (\bB) = \boldsymbol\mu_\alpha,   \text { for } \alpha \in F_\Gamma,
	\end{aligned}
\end{equation}
where the cost tensor decouples as in \cref{eq:cost_structure}, and $F_\Gamma \subset F$.
 Note that the graph structured MOT problem \cref{eq:ot_graph} is a special case of the standard MOT problem \cref{eq:omt_multi_discrete} in terms of the cost function. On the other hand, problem \cref{eq:ot_graph} allows for a more general class of constraints than the original problem \cref{eq:omt_multi_discrete}.

\begin{remark}
	We assume that $F_\Gamma \subset F$ for convenience of notation. However, in many cases there is not a cost associated with every factor in $F_\Gamma$. In this case we can simply define $C_\alpha(\bx_\alpha)=0$ for the corresponding $\alpha \in F_\Gamma$. In the examples and applications in this paper, we do not explicitly define these vanishing parts of the cost tensor. 
\end{remark}

The most widely studied graph-structured MOT problem is the Wasserstein  barycenter problem \cite{AguCar11}.
However, in a recent line of research, graph-structured MOT problems have been utilized in a wide range of applications.
For instance, in \cite{HaaRinChe20, haasler2021control} different control and estimation problems for large systems of indistinguishable agents are modeled as graph-structured MOT problems.
Even different classes of agents can be modeled with this framework, e.g., in network flow problems \cite{haasler2023scalable} and mean field games \cite{ringh2021graph, ringh2023mean}.
In the following we briefly review a few more exampled of optimal transport problems that fall into the category of graph-structured MOT.
In particular, we note that to the best of our knowledge, the Euler problem and Wasserstein splines have not yet been classified as such. 
Throughout this paper, we illustrate constrained marginals as gray nodes, and white nodes are estimated in the respective problems.

\subsubsection{Generalized Euler flow} \label{subsec:euler}
\vspace{-5pt}
One of the first applications of multi-marginal optimal transport, the generalized Euler flow \cite[Section~4.3]{BenCarCut15}, falls in fact into the category of a graph-structured problem.
The Euler flow studies a relaxation of the incompressible Euler equation, where one seeks trajectories of particles between constrained initial and final data, that minimize the kinetic energy \cite{brenier1989least}.
Discretizing the problem into $J$ time points, we describe the flow by a $J$-mode multi-marginal optimal transport plan $\bB \in \mR_+^{d^J}$, where each marginal describes the density of the flow at time instance $j$, discretized into $d$ points.
Note that the density of an incompressible flow, discretized uniformly over a given space, is described by a uniform distribution.
Thus, we constrain each marginal of the transport plan by $P_j (\bB) = \mu_j = \frac{1}{d} \boldsymbol{1}$.
In the Euler flow, it is enforced that each particle at time $0$ is transported to a given location at time $T$, given by a permutation $\sigma:\{1,\dots,d\} \to  \{1,\dots,d\}$.
This is modeled through a bi-marginal constraint
$P_{1,J} (\bB) = \boldsymbol{\mu}_{\{1,J\}} = \Pi_\sigma$, where $\Pi_\sigma \in \mR^{d\times d} $ is the permutation matrix for the permutation $\sigma$.
The Euler flow problem is thus formulated as the MOT problem \looseness=-1
\begin{equation} \label{eq:Euler1}
    \begin{aligned}
		\min_{\bB \in \mR_+^{d\times \dots \times d}} & \langle \bC, \bB \rangle                                                    \\
		\text{ subject to ~~}                               & P_j (\bB) = \mu_j   \text { ~~for } j =1,\dots,J; \quad P_{1,J} (\bB) = \boldsymbol{\mu}_{\{1,J\}},
	\end{aligned}
\end{equation}
where
\begin{equation} \label{eq:costEuler1}
    \bC (x_1,\dots,x_J) = \sum_{j=1}^{J-1} C (x_j,x_{j+1}),
\end{equation}
and the elements of the matrix $C \in \mR^{d\times d}$ contain the squared Euclidean distances between the discretization points, namely, $C_{x_j,x_{j+1}} = \|x_{j} - x_{j+1}\|^2$. 

\cite{BenCarCut15} also propose an alternative relaxed formulation, where the bi-marginal constraint in \cref{eq:Euler1} is replaced by an additional cost term in \cref{eq:costEuler1}, and the problem becomes
\begin{equation} \label{eq:Euler2}
    \begin{aligned}
		\min_{\bB \in \mR_+^{d\times \dots \times d}} \langle \bC, \bB \rangle                                                    
	\qquad	\text{ subject to }                               P_j (\bB) = \mu_j   \text { for } j =1,\dots,J,
	\end{aligned}
\end{equation}
where 
\begin{equation} \label{eq:costEuler2}
    \bC (x_1,\dots,x_J) = \sum_{j=1}^{J-1} C (x_j,x_{j+1}) + C^\sigma (x_1,x_J),
\end{equation}
and the elements of the matrix $C^\sigma \in \mR^{d\times d}$ are given by $C^\sigma (x_1,x_{J})  = \|x_{\sigma(1)} - x_{J}\|^2$.

\begin{figure}[tb] 
\caption{Illustration of graphical structures in some applications.  Gray nodes correspond to fixed marginals, and white nodes are estimated in the problem. The dotted edge between $x_1$ and $x_J$ in (a) may represent either a bi-marginal constraint, as in \cref{eq:Euler1}-\cref{eq:costEuler1}, or a cost interaction, as in \cref{eq:Euler2}-\cref{eq:costEuler2}.} \label{fig:examples}
    \vspace{-10pt}
    \begin{minipage}[t]{0.32\textwidth}
    \centering
\subcaption{Euler flow 
 }	\label{fig:euler_graph}
\begin{tikzpicture}
\scriptsize
  \tikzstyle{circ}=[circle, minimum size = 5mm, thick, draw =black!80, node distance = 3mm]
  \node[circ, fill=gray!50] (c1) {$x_1$};
  \node[circ, fill=gray!50] (c2) [right=of c1] {$x_2$};
  \node[circ, fill=gray!50,  node distance=8mm] (c4) [right=of c2] {$\!\!x_{J\!-\!1\!}\!\!$}; 
  \node[circ, fill=gray!50] (c5) [right=of c4] {$x_J$};
  \draw (c1) -- (c2);
  \draw[dotted] (c2) -- (c4);
  \draw (c4) -- (c5);
  \draw[densely dotted, thick] (c1) -- ++(0,1cm);
  \draw[densely dotted, thick] (c5) -- ++(0,1cm);
  \draw[densely dotted, thick] (c5) ++(0,1cm) -- ++(- 3.3,0cm);
\end{tikzpicture}
    \end{minipage}
    \begin{minipage}[t]{0.32\textwidth}
        \centering
	\subcaption{Wasserstein least squares
 }\label{fig:wasserstein_least_square}
\begin{tikzpicture}
\scriptsize
\tikzstyle{circ}=[circle, minimum size = 6mm, thick, draw =black!80, node distance = 3mm]
    \node[circ, fill=gray!50] (c1)  {$x_1$};
    \node[circ, fill=gray!50] (c2) [right=of c1] {$x_2$};
    \node[circ, fill=gray!50,  node distance=8mm] (cN) [right=of c2] {$x_J$}; 
    \node[circ] (c0) [above=4mm of c2, xshift=-6mm] {$x_0$}; 
    \node[circ] (cNp1) [above=4mm of c2, xshift=10mm] {$\!\!x_{\!J\!+\!1\!}\!\!$}; 
    \draw[loosely dotted] (c2) -- (cN);
    \draw (c0) -- (c1) -- (cNp1) -- (c0);
    \draw (c0) -- (c2) -- (cNp1);
    \draw (c0) -- (cN) -- (cNp1);
\end{tikzpicture}
    \end{minipage}
\begin{minipage}[t]{0.33\textwidth}
\centering
\subcaption{Wasserstein splines} \vspace{-3pt}\label{fig:wasserstein_spline}
    \begin{tikzpicture}
    \scriptsize
		\tikzstyle{circ}=[circle, minimum size = 6mm, thick, draw =black!80, node distance = 3mm]
		\node[circ, fill=gray!50] (c1) {$x_0$};
		\node[circ, fill=gray!50] (c2) [right=of c1] {$x_1$};
		\node[circ, fill=gray!50,  node distance=10mm] (c4) [right=of c2] {$\!\!x_{J\!-\!1\!}\!\!$};
		\node[circ, fill=gray!50] (c5) [right=of c4] {$x_J $};
		\draw (c1) -- (c2);
		\draw[dotted] (c2) -- (c4);
		\draw (c4) -- (c5);
		\node[circ] (c6) [above=of c1] {$v_0$};
		\node[circ] (c7) [above=of c2] {$v_1$};
		\node[circ] (c9) [above=of c4] {$\!\!v_{J\!-\!1}\!\!$};
		\node[circ] (c10) [above=of c5] {$v_J$};
		\draw (c6) -- (c7);
		\draw[dotted] (c7) -- (c9);
		\draw (c9) -- (c10);
		\draw (c6) -- (c2);
		\draw (c1) -- (c7);
		\draw (c4) -- (c10);
		\draw (c5) -- (c9);
		\draw (c1) -- (c6);
		\draw (c2) -- (c7);
		\draw (c4) -- (c9);
		\draw (c5) -- (c10);
	\end{tikzpicture}
\end{minipage}
    
\end{figure}
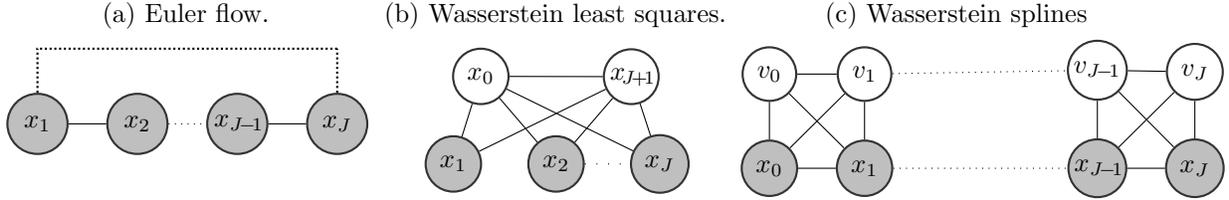

Both problems \cref{eq:Euler1}-\cref{eq:costEuler1} and \cref{eq:Euler2}-\cref{eq:costEuler2} can be described by a cyclic graph as illustrated in 
 \cref{fig:euler_graph}.
In particular, we note that although the authors in \cite{BenCarCut15} present only the computational method for the relaxed formulation \cref{eq:Euler2}-\cref{eq:costEuler2}, our proposed framework encompasses both formulations. 

\subsubsection{Wasserstein least squares}\label{sec:wls}
\vspace{-5pt}
Consider a least square fitting problem in the Wasserstein-2 space. Given a set of probability measures $\rho_{t_1}, \rho_{t_2},\ldots,\rho_{t_J}$ with $0\le t_1<t_2 < \cdots< t_J \le 1$, we aim to find two distributions $\mu_0$ and $\mu_1$ such that the set of measures is closest to the displacement interpolation between $\mu_0$ and $\mu_1$ in terms of Wasserstein-2 distance. More specifically, we seek the solution to
\begin{equation}\label{eq:W2leastsquare}
	\min_{\mu_0,\mu_1} \sum_{j=1}^J W_2 (\rho_{t_j}, \mu_{t_j })^2
\end{equation}
where $\mu_t$ denotes the displacement interpolation connecting $\mu_0$ and $\mu_1$. This least square problem can be reformulated as
\begin{equation}
\begin{aligned}
\min_{\pi, \mu_0, \mu_1} & \int c(x_0,x_1,\cdots,x_J , x_{J+1}) \d\pi (x_0,x_1,\cdots,x_J,x_{J+1}) \label{eq:W2leastsquare_multimarginal}
\\ 
\text{subject to } & 
P_j (\pi) = \rho_{t_j },~~j=1,2,\cdots,J;  \qquad 
P_0 (\pi) = \mu_0 , \qquad 
P_{J+1} (\pi) = \mu_1 . 
\end{aligned}
\end{equation}
By discretizing the integral in \cref{eq:W2leastsquare_multimarginal} to a summation, we can transform the problem into the form presented in \cref{eq:ot_graph}.
To enforce that the optimality of the coupling between $\mu_0$ and $\mu_1$, we can add a large penalty on it, leading to the modified cost
	\begin{equation}\label{eq:wls_cost}
		c(x_0,x_1,\cdots,x_J,x_{J+1}) = \sum_{j=1}^J \|x_j-(1-t_j)x_0-t_i x_{J + 1}\|^2 +\alpha \|x_0-x_{J +1}\|^2
	\end{equation}
	with $\alpha>0$ being a sufficiently large number. 


\subsubsection{Wasserstein Splines} \label{subsec:splines}
\vspace{-5pt}
The Wasserstein spline developed in \cite{CheConGeo18} is a natural notion of cubic spline over the space of probability distributions.
Given a set of points $x_{t_0}, x_{t_2}, \ldots, x_{t_J}$ with $0= t_0<t_1<\cdots < t_J=1$, a cubic spline is a piece-wise cubic polynomial $\{x(t)\,\mid\, 0\le t\le 1\}$ that is smooth and matches the given data points in the sense that $x(t_j) = x_{t_j}$ for all $j = 0,1,\ldots,J $. The cubic spline can also be characterized as the interpolating curve of the data points with minimum mean square acceleration. Based on this characterization, \cite{CheConGeo18} generalized the notion of cubic spline to interpolate a set of probability distributions $\rho_{t_0}, \rho_{t_1}, \ldots, \rho_{t_J}$. The definition has a static reformulation in the phase space as
\begin{equation*}
\begin{aligned}
\min_\pi  \int c(x_0,v_0,\cdots,x_J,v_J) \d\pi(x_0,v_0,\cdots,x_J,v_J) \quad \text{ s.t. } 
P_{x_j}(\pi) = \rho_{t_j},\ j = 0,1,\ldots, J,
\end{aligned}
\end{equation*}
where $c(x_0,v_0,\cdots,x_J,v_J) = $
\begin{equation*}
 \sum_{i=0}^{J-1}\frac{1}{t_{\!j + 1\!}-t_j} \!  \left( \! 12 \left\|\frac{x_{\!j + 1\!}-x_j}{t_{\!j + 1\!}-t_j}-v_j \right\|^2 \!\! -12 \left\langle \frac{x_{\!j + 1\!}-x_j}{t_{\!j + 1\!}-t_j}-v_j,v_{\!j + 1\!}-v_j \right\rangle + 4\|v_{\!j + 1\!}-v_j\|^2 \! \right) \!\!. 
\end{equation*}
This problem is associated with the graphical model in \cref{fig:wasserstein_spline}, where the constraints are on the position nodes $x_0,x_1,\ldots,x_J$.


\subsection{Multi-marginal optimal transport with global regulatization} \label{sec:global}

One approach for solving multi-marginal optimal transport problems with a graph-structured cost \cref{eq:cost_structure} builds on regularizing the problem with a global entropic regularization term as in \cref{eq:omt_multi_regularized} \cite{HaaSinZha20,fan2021complexity,HaaRinChe20,SinHaaZha20,ringh2023mean}.
In this section we briefly present the core idea of this approach, since the parallel algorithms that we develop in the following sections builds on similar concepts.
For a deeper introduction to the topic we refer the reader to the supplementary material, where we present the framework concisely, and in a more general setting than previous works.

It turns out that in case the underlying graph is a tree, i.e., it does not have any loops, computing a projection $P_j(\bK \odot \bU \odot \bM)$, as needed in Algorithm \ref{alg:sinkhorn}, scales linearly in the number of marginals \cite{HaaSinZha20}, whereas a brute force approach for computing the projections scales exponentially.
This result can be extended to general graphical models with possible loops, by partitioning the graph as a so-called junction tree \cite{fan2021complexity}. A junction tree (also called tree decomposition) describes a partitioning of a graph, where several nodes are clustered
together, such that the interactions between the clusters can be described by a tree structure.

\begin{definition}[Junction tree \cite{dawid1992applications}] \label{def:junction_tree}
	A junction tree $\mathcal{T}= (\mathcal{C}, \mathcal{E})$ over a graph 
 {$G=(V,F)$}
 is a tree whose 
 nodes are cliques
 $c\in \mathcal{C}$ associated with subsets $\bx_c \subset V$, and 
 satisfying
 the 
 properties:
	\begin{itemize}
  \item Family preservation: For each $\alpha \in F$ there is a clique $ c \in \cC$ such that 
  $\alpha \subset c $.
\item Running intersection: For every pair of cliques $c_i, c_j \in \mathcal{C}$, every clique on the path between $c_i$ and $c_j$ contains $c_i \cap c_j$. 
	\end{itemize}
	For two adjoining cliques $c_i$ and $c_j$, define the separation set $  s_{ij}= \{ v\in V : v\in c_i \cap c_j \} $. 
\end{definition}

 A graph typically has several possible junction trees. In particular, one trivial junction tree is always given by a single cluster containing all nodes.
For the computational method in \cite{fan2021complexity} the junction tree  is set up such that each constraint marginal or clique is represented by a leaf in the junction tree.
The entropy-regularized graph structured optimal transport problem with junction tree decomposition $\mathcal{T}= (\mathcal{C}, \mathcal{E})$ is thus of the form 
\begin{equation} \label{eq:ot_jt}
		\min_{\bB \in \mR_+^{ d\times \dots \times d
  }} \langle \bC, \bB \rangle      + \epsilon \cH(\bB 
  {| \bM}  ) 
	\qquad	\text{ subject to }                               P_c (\bB) = \boldsymbol\mu_c,  \text { for } c \in \mathcal{C}_\Gamma,
\end{equation}
where $\mathcal{C}_\Gamma \subset \mathcal{C}$ is the set of constrained cliques,
and the cost tensor decouples as
\begin{equation}\label{eq:cost_jt}
	\bC(\bx) = \sum_{c \in \mathcal{C}} \bC_{c}(\bx_c).
\end{equation}
Such a junction tree representation exists for any graph, however it may require adding artificial clusters to the junction tree decomposition.
Junction trees for the examples in Sections \ref{subsec:euler}--\ref{subsec:splines} are sketched in the supplementary material.

\section{
	Tree-structured multi-marginal optimal transport with local regularization}\label{sec:tree_local_reg}


In the context where the MOT problem is tree structured, i.e., structured as a graph without cycles, the innate characteristics of tensor inner products facilitate a natural reformulation of the problem. This reformulation manifests as a sum of coupled bi-marginal OT problems. Considering a set of fixed marginals, denoted as $\mu_j$ for each $j$ belonging to the vertex set $\Gamma$, the problem of coupled bi-marginal OT can be written as: 
\begin{equation} \label{eq:ot_graph_pairwise}
    \min_{\substack{ B_{(j,k)}, (j,k)\in E                                                      \\ \mu_j, j\in V\setminus \Gamma}}  \sum_{(j,k)\in E} 
\<C_{(j,k)} , B_{(j,k)} \>
     \quad   \text{ s.t. }  B_{(j,k)} \mathbf{1} = \mu_{j} , 
        ~~ B_{(j,k)}^\top \mathbf{1} = \mu_{k} ~~ \text{for } (j,k)\in E ,
\end{equation}
where $E$ is a set of all the edges of the graph.
The matrices $C_{(j,k)} \in \mathbb{R}^{d \times d}$ and $B_{(j,k)} \in \mathbb{R}^{d \times d}$ represent the transport cost and transport plan corresponding to the edge connecting nodes $j$ and $k$, respectively. 
We implicitly assume the convention $B_{(j,k)} =B_{(k,j)}^\top$ and $C_{(j,k)} =C_{(k,j)}^\top$.
In contrast to the global entropy regularization presented in \cref{sec:global} for a standard MOT problem \cref{eq:ot_jt}, applying local regularization on each edge appears as a more direct and natural
approach within this framework. Notably, in tree-structured optimal transport scenarios, such as barycenter problems, this form of regularization is frequently employed, as described in existing literature \cite{kroshnin2019complexity}.


\subsection{Algorithm}\label{sec:tree_algo}

Drawing inspiration from \cite{kroshnin2019complexity}, we developed a {\em parallel} algorithm to address the coupled bi-marginal OT problem for general trees. 
Herein, we propose to use a
modified entropy regularization 
similar to \cite{fan2021complexity}. The reference joint distribution between nodes $j$ and $k$ is denoted as:
\begin{equation} \label{eq:jointdist}
M_{(j,k)} := \gamma_j \otimes \gamma_k,
\end{equation}
where $\gamma_j = \mu_j $ if $j \in \Gamma$ and 
$\mathbf{1}$ otherwise. 
As discussed before, we always assume that marginals may only be constrained if they are on the leaves of the tree.
It is important to note that therefore the scenario where $(j,k) \in E$ and both $j, k \in \Gamma$ does not exist. In such cases, the OT problem simplifies into several independent subproblems. Consequently, this leads us to formulate the regularized problem with an enforced mass constraint. 
\begin{equation} \label{eq:ot_graph_pairwise_mass}
    \begin{aligned}
        \min_{\substack{ B_{(j,k)}, (j,k)\in E \\ \mu_j, j\in V\setminus \Gamma}} \quad & \sum_{(j,k)\in E} \langle C_{(j,k)}, B_{(j,k)}\rangle + \epsilon \cH(B_{(j,k)} |M_{(j,k)}  ) \\
        \text{subject to} \quad & B_{(j,k)} \mathbf{1} = \mu_{j}, ~ 
          B_{(j,k)}^\top \mathbf{1} = \mu_{k}  \text{ for } (j,k)\in E; 
          ~~~  \mathbf{1}^\top \mu_j = 1 
           \text{ for } j \in V \backslash \Gamma.
    \end{aligned}
\end{equation}
The regularization, as formulated, is applied locally to each edge of the tree, 
in contrast to the previous formulations \cref{eq:ot_jt}, which regularizes the multi-marginal optimal transport tensor.
Utilizing a standard Lagrangian dual theory enables the derivation of the dual problem.
Denote 
\begin{equation}\label{eq:bi-marginal_B_sol}
 B_{(j,k)} = \diag (u_{(j,k)} \odot \gamma_j ) K_{(j,k)} \diag (u_{(k,j)} \odot \gamma_k ), \qquad K_{(j,k)} = \exp(- C_{(j,k)}/ \epsilon)
\end{equation}
where $u_{(j,k)} \in \mR^{d} $ .
Furthermore, $u=\{u_{(j,k)} \}_{j \neq k}$ and $\rho= \{\rho_j \}_{j \in V\backslash \Gamma}$. For a node $j \in \Gamma$, its unique neighbor is denoted by $k_j$. The dual problem is then formulated as:
\begin{equation} \label{eq:tree_dual_mass_consv}
			\begin{aligned}
\min_{u,  \rho } ~~f(u , \rho ):= &  \sum_{(j,k)\in E} \|B_{(j,k)}\|_1 - \sum_{j \in \Gamma} \log (u_{(j,k_j)})^\top \mu_{j} - \sum_{j \in V \backslash \Gamma} \rho_j  \\
\text{subject to }                                 & \ \sum_{k \in N(j)} \log(u_{(j,k)}) = \rho_j \mathbf{1}, \quad \text{for } j \in V \backslash \Gamma.
			\end{aligned}
		\end{equation}
We then introduce another set of dual variables 
and apply the Lagrangian dual theory again 
to the dual problem \cref{eq:tree_dual_mass_consv}, and obtain the update rules. The dual variables we introduce are reflected in $v_j$ in \cref{eq:free_node_update_mass_cons}.

\begin{theorem}\label{thm:tree_dual_var}
The optimal dual variables can be found as the limit point of the following iterative scheme:
For $j\in \Gamma$ we update
\begin{equation}\label{eq:fixed_node_update_mass_cons}
        u_{(j,k_j)} = \mathbf{1} ./ (K_{(j,k_j)} u_{(k_j,j)}).
    \end{equation}
For $j\in V \setminus \Gamma$ we update the set of vectors $u_{(j,k)}$, where $k\in N(j)$, according to 
\begin{equation}
\label{eq:free_node_update_mass_cons}
			u_{(j,k)} = e^{\rho_j /|N(j)|}  v_j ./ (K_{(j,k)} (u_{(k,j)} \odot \gamma_k)),
  \quad
    v_j \! = \bigg( \bigodot_{k\in N(j)} \! \Big( K_{(j,k) } (u_{(k,j)} 
   \odot \gamma_k ) \Big)
 \bigg)^{\!1/|N(j)|} \!
		\end{equation}
and we update
$\rho_j = - |N(j)| \log \left(\mathbf{1}^\top v_j \right).$
\end{theorem}

Based on \cref{thm:tree_dual_var}, we note that for any node $j$, the updates of $\{u_{(j,k)}\}_{k \in N(j)} $ or $\rho_j$ depend only on the variables $u_{(k, j)}$ from neighboring nodes $k \in N(j)$. This observation is critical for enabling parallel updates in the algorithm, as it allows for independently updating dual variables associated with two distinct parts of the tree. The concept of a bipartite graph becomes instrumental here, providing the framework for such parallelism.

\begin{figure}[htb]
\noindent 
\begin{minipage}[t]{0.5\linewidth}
    \begin{definition}{(Bipartite Graph \cite{metcalf2016cybersecurity})}.  
    A bipartite graph is a graph where the vertices can be divided into two disjoint sets such that all edges connect a vertex in one set to a vertex in another set. There are no edges between vertices in the disjoint sets.
    \end{definition}
\end{minipage} \hfill
\begin{minipage}[t]{0.46\linewidth}
    \centering
    \vspace{-10pt}
    \captionof{figure}{Bipartite graph with sets $S_1,$ $S_2$.}
    \label{fig:bipartite_graph}
    \vspace{-10pt}
    \scalebox{0.8}{
    \begin{tikzpicture}[node distance=1cm and 0.5cm]    
        \tikzstyle{circ}=[circle, minimum size = 5mm, thick, draw=black!80, fill=white, node distance = 12mm]
        \tikzstyle{block}=[rectangle, draw=black!80, dashed, inner sep=2mm]

        \foreach \x in {1,...,4} {
            \node[circ] (a\x) at (\x,0) {};
        }
        \foreach \x in {1,...,3} {
            \node[circ] (b\x) at (\x+0.5,-1.5) {};
        }
        \foreach \i in {1,...,4} {
            \foreach \j in {1,...,3} {
                \draw (a\i) -- (b\j);
            }
        }
        \node[block, fit=(a1) (a4), label={[label distance=2mm]left:$S_1$}] {};
        \node[block, fit=(b1) (b3), label={[label distance=2mm]left:$S_2$}] {};
    \end{tikzpicture} 
    }
\end{minipage}
\end{figure}

To leverage this, the algorithm unfolds in the following manner:
We begin by representing the tree as a bipartite graph, a natural representation since all trees are inherently bipartite according to \cite{scheinerman2012mathematics}.
The algorithm then iteratively updates until convergence, alternating between:
i) Updating all dual variables from the first partition to the second.
ii) Updating all dual variables from the second partition to the first.
Denote the intermediate computational variables
\begin{equation} \label{eq:q}
 q^{t}_{(j,k)} = B^{t}_{(j,k)} \mathbf{1} , \qquad   q^{t+1}_j = \bigg(\bigodot_{k \in N(j)} q^{t}_{(j,k)}\bigg)^{\frac{1}{|N(j)|}} .
\end{equation}
Our method is detailed in \cref{alg:tree_local_reg}, termed \texttt{Bipartite Iterative Scaling}. 
A supporting \cref{cor:equi} underpins its theoretical foundation.

\begin{algorithm}[tb]
\KwIn{ Probability vectors $\{\mu_j\}_{j \in \Gamma}$. Stopping criteria $\delta'$. 
Transport cost $ \{ C_{(j,k)} \}
$. Regularization 
param.
$\epsilon$.
}

{Construct a bipartite graph where the nodes are divided into two disjoint groups $\{S_1, S_2\}$.}	

{Initialization of variables:} $t=1, ~~u^0_{(j,k)} = \mathbf{1} \in \mR^d $ for $ j \neq k$

\While{$\sum_{j \in \Gamma} \|B_{(j,k_j)}^t \mathbf{1} - \mu_j\|_1  +  \sum_{j \notin \Gamma}  \sum_{k \in N(j)} \|q_{(j,k)}^t - \bar q_j^t\|_1 \ge \delta' $ }{ 
\vspace{0.2cm}
Let $S =S_1$  if $t$ is odd, otherwise $S =S_2$

Calculate $q_{\{j,k\}}^t, ~~ q_j^{t+1}$ according to \cref{eq:q}  for $j \in S $
\vspace{-5pt}
\begin{equation} 
u^{t+1}_{(j,k_j)}  = u^{t}_{(j,k_j)} \odot \mu_j ./ q^t_{(j,k_j)}
\qquad \text{for } j \in \Gamma \bigcap S \label{eq:bipartite_tree_updates1}
\end{equation}
\vspace{-10pt}
\begin{equation} \label{eq:bipartite_tree_updates2} 
u_{(j,k)}^{t+1} = u_{(j,k)}^t 
\odot \frac{q_j^{t+1}}{\|q_j^{t+1}\|_1}  ./  q_{(j,k)}^t  \qquad \text{for } j \in (V \backslash \Gamma) \bigcap S 
\end{equation}
\vspace{-10pt}

$t  \leftarrow t+1$

Construct ${B}^t_{(j,k)}, ~ (j,k) \in E $ according to \cref{eq:bi-marginal_B_sol}.
  }

\Return{${B}^t_{(j,k)}$ for all 
 $ (j,k) \in E $ }
\caption{
Bipartite Iterative Scaling Algorithm
for graphical coupled bi-marginal
OT 
with local regularization
	} 	
\label{alg:tree_local_reg}
\end{algorithm}

\begin{algorithm}[tb]
\tcc{
Post-update of nodes in one partition, the projections on nodes in the opposite partition do not match the marginals. This partition is denoted as $S$.
}	
	\KwIn{The transportation plans ${B}_{(j,k)}$ for all $(j,k)\in E$, where $P_{(j,k)}$ has mismatch with marginal $q_j$ when $j \in S$. A sequence of probability vectors $\{\mu_j\}_{j \in \Gamma}$
	}
\For{$j \in S$}{ 
  \If{$j \in \Gamma$}{Round  $B_{(j, k_j)}$ by \cite[Algorithm 2]{altschuler2017near} such that $\widehat{B}_{(j, k_j)} \in \Pi(\mu_j, B^\top_{(j, k_j)} \ett ) $ . }
  \Else{Round  $B_{(j, k)}$ by  \cite[Algorithm 2]{altschuler2017near} such that $\widehat{B}_{(j, k)} \in \Pi(\bar q_j(B) ,B^\top_{(j, k)} \ett ) $ for $\forall k \in N(j)$, where  $\bar q_j (B) = \frac{1}{|N(j)|} \sum_{k \in N(j)} B_{(j,k)}  \mathbf{1}$.
  }
  }
\Return{$\widehat{B}_{(j,k)}$ for all $(j,k)\in E$}
	\caption{Rounding for 
  transport plan on a bipartite graph
 } 	
\label{alg:round}
\end{algorithm}

\begin{proposition} \label{cor:equi}
The updates in \cref{thm:tree_dual_var} can be equivalently expressed as the updates \cref{eq:bipartite_tree_updates1}-\cref{eq:bipartite_tree_updates2} in \cref{alg:tree_local_reg} 
together with the variables \cref{eq:q}. 
\end{proposition}

The updates of the dual variables, as specified in \cref{eq:bipartite_tree_updates1}-\cref{eq:bipartite_tree_updates2}, are independent within the same group $S$, facilitating simultaneous updates in {\em parallel}. Provided there is sufficient memory to store all variables, concurrent updating is feasible.

In each iteration \cref{alg:tree_local_reg} updates the variables corresponding to one of the partitions $S_1$ or $S_2$.
This update guarantees that the constraints are matched exactly for the nodes in this partition. 
However, this means that after any finite number of iterations of \cref{alg:tree_local_reg} the marginal constraints will not be matched for the nodes in the partition that has not been updated in the last iteration. 
To address this, we adopt a strategy akin to the ones in \cite{altschuler2017near, fan2021complexity} and employ a Rounding \cref{alg:round}.  This ensures the generation of transport plans that adhere to all stipulated constraints in \cref{eq:ot_graph_pairwise_mass}.
 The complexity of the rounding process is $ \cO( |E|d^2  )$. 
 \cref{alg:round} differs from the method in \cite{fan2021complexity} in that it requires modification of the transport plans across all edges.


\subsection{Complexity analysis}\label{sec:tree_complexity}

In this section, we provide a rigorous demonstration that our \cref{alg:tree_local_reg}, when used in conjunction with the Rounding \cref{alg:round},  requires $T = \widetilde{\mathcal{O}} \left( {|E|^2}/{\delta^2} \right)$ 
iterations  to achieve a $\delta$-approximation of the problem \cref{eq:ot_graph_pairwise}. First, we will illustrate that, through the specific mass normalization technique employed in \cref{eq:bipartite_tree_updates2}, the transport plan maintains unit mass in each iteration.


\begin{lemma}\label{cor:tree_mass}
At each update $t$ of \cref{alg:tree_local_reg}, we have
$\|B^t_{(j,k)}\|_1 =  \|q_{(j,k)}^t\|_1 = 1, ~~\forall  (j,k) \in E .$
\end{lemma}

We now define the variable $\lambda^t_{(j,k)} =\epsilon \log (u^t_{(j,k)}) $ for $(j,k) \in E$, and $u^t_{(j,k)}$ are the iterates in Algorithm \ref{alg:tree_local_reg}.
\cref{lem:tree_dual_varia_bounds} provides bounds on the scaling variables.
\begin{lemma} \label{lem:tree_dual_varia_bounds}
Denote $C_\infty : =\max_{(j,k) \in E }\|C_{(j,k)}\|_\infty.
$
For all iterates $\lambda_{(j,k)}^t$,
and their fixed points $\lambda_{(j,k)}^*$,
\begin{equation*}
 \max_{i=1,\dots,d} (\lambda_{(j,k)}^t)_i - \min_{i=1,\dots,d} (\lambda_{(j,k)}^t)_i  \leq2 C_\infty , \qquad  \max_{i=1,\dots,d} (\lambda_{(j,k)}^*)_i - \min_{i=1,\dots,d} (\lambda_{(j,k)}^*)_i  \leq2 C_\infty .
\end{equation*}
\end{lemma}



The following lemma links the dual objective error in \cref{thm:tree_dual_var} to  the stopping criteria in \cref{alg:tree_local_reg}.

\begin{lemma}\label{lem:f_upper_bound_tree}
For all iterates $(u^t,\rho^t)$ of the method in \cref{alg:tree_local_reg} and their fixed points $(u^*,\rho^*)$, it holds that 
\begin{equation*}
    f(u^t,\rho^t) - f(u^*,\rho^*)
    \le 
    \frac{2 C_\infty}{\epsilon} \left( \sum_{j \in {\Gamma}} \|B_{(j,k_j)}^t \mathbf{1} - \mu_j\|_1 
    +\sum_{j \in {V\backslash \Gamma}} \sum_{k \in N(j)} \| B_{(j,k)}^t \mathbf{1} - \bar q_j^t\|_1 \right)
\end{equation*}
where $\bar q_j^t= \frac{1}{|N(j)|} \sum_{k \in N(j)} q_{(j,k)}^t= \frac{1}{|N(j)|} \sum_{k \in N(j)} B_{(j,k)}^t  \mathbf{1}$.
\end{lemma}
The increase in the dual objective between iterations is also closely related to  the stopping criteria in \cref{alg:tree_local_reg}.

\begin{lemma}\label{lem:f_lower_bound_tree} 
For the iterates in \cref{alg:tree_local_reg}, it holds that
\begin{equation*}
 f(u^t,\rho^t) - f(u^{t+1} \!, \rho^{t+1}) 
    \ge  \frac{1}{22|E|} \! \left( \sum_{j \in {\Gamma}} \|B_{(j,k_j)}^t \mathbf{1} - \mu_j\|_1 
    + \!\! \sum_{j \in {V\backslash \Gamma}} \sum_{k \in N(j)} \!\! \| B_{(j,k)}^t \mathbf{1} - \bar q_j^t\|_1 \!\right)^{\!\!2} \!\!.
\end{equation*}
\end{lemma}

Unlike \cite{kroshnin2019complexity}, our \cref{lem:f_upper_bound_tree} and \cref{lem:f_lower_bound_tree} apply universally for any $t$, obviating the need to consider odd and even $t$ distinctly.

\begin{theorem}\label{thm:num_iter}
\Cref{alg:tree_local_reg} generates $\{B_{(j,k)}^t\}$ satisfying 
$\sum_{j \in \Gamma} \|B_{(j,k_j)}^t \mathbf{1} - \mu_j\|_1  +  \sum_{j \notin \Gamma}  \sum_{k \in N(j)} \|q_{(j,k)}^t - \bar q_j^t\|_1   \le \delta'$
within 
    $t  = 2+ \frac{88 |E| C_\infty}{\delta' \epsilon}$
iterations. 
\end{theorem}

The subsequent two lemmas quantify the error introduced to the objective value by the rounding process.

\begin{lemma}\label{lem:round}
Let $\{B_{(j,k)}\}$ be a sequence of transportation plans on a bipartite tree. 
\cref{alg:round} returns $\{\widehat B_{(j,k)}\}$ satisfying $\widehat B_{(j,k)}^\top \ett= q_j $ for $ j \in V $, {where} $q_j$ is the marginal {distribution} of node $j$.
Moreover, it holds that
\begin{equation*}
    \begin{aligned}
& \bigg|  \sum_{(j,k) \in E } \< C_{(j,k)}, B_{(j,k)} \> 
- \sum_{(j,k) \in E } \< C_{(j,k)}, \widehat B_{(j,k)} \> \bigg| \\
& \le  2 C_\infty  \left( \sum_{j \in \Gamma }   \| \mu_j - B_{(j,k_j)} \ett) \|_1  + \sum_{j \notin \Gamma } \sum_{k \in N(j) }   \| \bar q_j (B) - B_{(j,k)} \ett \|_1 \right) ,
    \end{aligned}
\end{equation*}
where
$\bar q_j (B) = \frac{1}{|N(j)|} \sum_{k \in N(j)} B_{(j,k)}  \mathbf{1}$.
\end{lemma}

\begin{lemma}\label{lem:round_applied}
Let $\{\widetilde B_{(j,k)}\}$  be the output of \cref{alg:tree_local_reg}
and let $\{\widehat B_{(j,k)}\}$ be the output of \cref{alg:round} with input $\{\widetilde B_{(j,k)}\}$, and denote $\{ B^*_{(j,k)}\}$ the optimal solution to the unregularized MOT problem \cref{eq:ot_graph_pairwise}. Then it holds that
\begin{equation*}
\begin{aligned}
& \sum_{(j,k) \in E } \< C_{(j,k)}, \widehat B_{(j,k)} \> - \sum_{(j,k) \in E } \< C_{(j,k)},  B^*_{(j,k)} \> \\
\le & 2 \epsilon |E|  \log(d) + 4 C_\infty  \left( \sum_{j \in \Gamma }   \| \mu_j - \widetilde B_{(j,k_j)} \ett) \|_1  + \sum_{j \notin \Gamma } \sum_{k \in N(j) }   \| \bar q_j (\widetilde B) - \widetilde B_{(j,k)} \ett \|_1 \right).
\end{aligned}
\end{equation*}
\end{lemma}

Given a predetermined final error threshold $\delta$, we select parameters $\epsilon = \frac{\delta}{4 |E| \log(d)}$ and $\delta' = \frac{\delta}{8 C_\infty}$ for use in \cref{alg:tree_local_reg}. This results in an output of $\widetilde{B}_{(j,k)}$, for each edge pair $(j,k) \in E$. Next, we employ \cref{alg:round} to round $\widetilde{B}_{(j,k)}$ for each $(j,k) \in E$, thereby producing $\widehat{B}_{(j,k)}$.
With these procedures in place, we are equipped to compute the final complexity for the coupled bi-marginal OT problem as described in \cref{eq:ot_graph_pairwise}. \looseness=-1


 
\begin{theorem} \label{thm:tree_local_reg_complexity}
\cref{alg:tree_local_reg}, in conjunction with \cref{alg:round}, yields a $\delta$-approximate solution to the tree-structured MOT problem \cref{eq:ot_graph_pairwise} in $T$ arithmetic operations, where
     $T = \mathcal O \left( \frac{|E|^3  C_\infty^2 d^2 \log(d)}{\delta^2 } \right) .$
\end{theorem}


\begin{remark}[{Comparison of complexity with global regularization} 
\cite{fan2021complexity}]\label{rem:complexity}
We compare our complexity results with the ones in \cite{fan2021complexity}, for graph-structured MOT with global entropic regularization (cf. supplementary material) in terms of the total complexity and the number of iterations.
i) \emph{Total Complexity}: According to Theorem 2 in \cite{fan2021complexity}, global entropy regularization requires $\widetilde \cO( d(G) |E| |\Gamma|^2 C_\infty^2 d^2 / \delta^2)$ arithmetic operations with high probability, where $d(G)$ denotes the maximum node distance in graph $G$. Under the assumptions $|E| = \cO( |\Gamma| )$ and $d(G) = \cO(1)$, akin to the Wasserstein barycenter problem, both the global and local regularization methods offer similar complexity bounds. However, a notable distinction is that \cref{thm:tree_local_reg_complexity} offers a deterministic bound.
ii) \emph{Iterations}: \cref{alg:tree_local_reg} with local regularization, as per \cref{thm:num_iter} and \cref{thm:tree_local_reg_complexity}, necessitates $\mathcal O \left( \frac{|E|^2 C_\infty^2 \log(d)}{\delta^2 } \right)$ iterations, less than the requirement for global regularization at $\cO \left( \frac{|E| |\Gamma|^2 C_\infty^2 \log(d)}{\delta^2 } \right)$. The reduced iteration count is a benefit of algorithm parallelization, albeit with a higher computational load per iteration.
\end{remark}


\section{
Generalized algorithm with local regularization to general graphs} \label{sec:extension}


Building on the concepts from \cref{sec:tree_local_reg}, we extend our approach to OT problems with a general factor graph structure as indicated in \cref{eq:cost_structure}. We employ a strategy that breaks down the 
MOT problem into interconnected smaller MOT problems, introducing local regularization for each sub-problem. To facilitate this, we define a specialized type of junction tree, elaborated in \cref{def:junction_tree}.






\begin{definition} \label{def:jt_general}
    A {\textbf{modified}} junction tree is a class of junction tree 
    $\mathcal{T}= (\mathcal{C}, \mathcal{E})$ 
    with two disjoint sets of cliques $ \cQ_C \subset \cC$ and separators $\cQ_S \subset \cC$ such that
    \begin{enumerate}[noitemsep, leftmargin=*]
        \item 
        For any $c \in \cQ_S$, its neighbors $N(c) \subseteq \cQ_C $, and for any  $c \in \cQ_C$, its neighbors $N(c) \subseteq \cQ_S $.
        \item Family preservation: For each $\alpha \in F$ there is a clique $ c \in \cQ_C$ such that 
        {$\alpha \subset {c}$.}
        \item  Running intersection: For every pair of cliques $c_i, c_j \in \cQ_C$, every clique in $\cQ_C$ or $\cQ_S$
        on the path between $c_i$ and $c_j$ contains 
        {$c_i \cap c_j$.}
        \item One clique is only connected to two separators: For all $c \in \cQ_C$ it holds $|N(c)|=2$.
    \end{enumerate}
\end{definition}
{Note that each separator set in \cref{def:junction_tree} is connected to exactly two cliques. On the other hand, by \cref{def:jt_general}, our cliques in $\cQ_S$ can be connected to an arbitrary number of cliques, and instead, the cost-cliques in $\cQ_C$ can only be connected to two cliques. This latter assumption is necessary for efficient parallelization of the algorithm we introduce and can be relaxed as we will explain in \cref{rem:more_neighbours}. }



\begin{example}\label{ex:mjt}
The junction trees in Figures~\ref{fig:traffic_modified_jt}, \ref{fig:wls_modified_jt}, and \ref{fig:spline_modified_jt} are \textbf{modified} junction trees for the respective graphs. In contrast to junction trees for the problem with global regularization, depicted in Figures~\ref{fig:traffic_graph_jt}, \ref{fig:W2leastsquare_jt}, and \ref{fig:spline_jt}, each clique (shown as circles) can only be connected to two separators (shown as squares), while separator can connect to multiple cliques.
Moreover, note that constraints in the respective MOT problems are indicated by gray boxes in Figures~\ref{fig:traffic_modified_jt}, \ref{fig:wls_modified_jt}, and \ref{fig:spline_modified_jt}. This is explained in more detail in the following section.
\end{example}

\begin{figure}[tb]
    \centering
      \caption{Modified junction tree for the graph of Euler flow example from \cref{fig:euler_graph}.
    }
    \label{fig:traffic_modified_jt}
\scalebox{0.7}{
    \begin{tikzpicture}
        \tikzset{
            circ/.style={circle, minimum size = 12mm, thick, draw=black!80, node distance=20mm},
            rect/.style={rectangle, draw=black!100, minimum width=12mm, minimum height=10mm},
            srect/.style={rectangle, draw=black!100, fill=gray!50, minimum width=8mm, minimum height=8mm},
            nrect/.style={rectangle, draw=black!100,  minimum width=8mm, minimum height=8mm},    
            inner/.style={rectangle, draw=black!100, fill=gray!50, minimum size=5mm, inner sep=0.5pt}
        }
        
        \node[srect, left=0.8cm of c1] (r0) {$2$}; 
        \node[circ] (c1) {$1,2,3$};
        \node[rect, right=0.5cm of c1] (r1) {}; 
        \node[inner] at ([xshift=2mm] r1.center) (i1) {$3$}; 
        \node[left=1mm of i1] {$1$}; 
        \node[circ, right=0.5cm of r1] (c2) {$1,3,4$};
        \node[rect, right=0.5cm of c2] (r2) {}; 
        \node[inner] at ([xshift=2mm] r2.center) (i2) {$4$}; 
        \node[left=1mm of i2] {$1$}; 
        \node[circ, right=0.5cm of r2] (c3) {$1,4,5$};
        \node[circ, right=0.5cm of c3] (c4) {$1,J\!-\!1,J$};
        \node[nrect, right=0.5cm of c4] (r3) {$1,J$}; 
        
        \draw (c1) -- (r0);
        \draw (c1) -- (r1);
        \draw (r1) -- (c2);
        \draw (c2) -- (r2);
        \draw (r2) -- (c3);
        \draw[dotted] (c3) -- (c4);
        \draw (c4) -- (r3);
    \end{tikzpicture}
    }
\end{figure}

\begin{figure}[tb]
\centering
\caption{Modified junction tree for the graph of Wasserstein least square example from \cref{fig:wasserstein_least_square}.}
\label{fig:wls_modified_jt}
\scalebox{0.7}{
\begin{tikzpicture}
    \tikzstyle{circ}=[circle, minimum size = 14mm, thick, draw =black!80, node distance = 6mm]
    \tikzstyle{crect}=[rectangle, minimum size = 6mm, thick, draw =black!80, node distance = 6mm] 
    \tikzstyle{rect}=[rectangle, minimum size = 6mm, thick, draw =black!80, fill=gray!50, node distance = 6mm] 
    \node[circ] (c1) {\Large $\substack{0, \ 1\\ J+1}$};
    \node[circ] (c2) [right=of c1]{\Large $\substack{0,\ 2\\ J+1}$};
    \node[circ] (c3) [right=of c2] {\Large $\substack{0,\ 3\\ J+1}$};
    \node[] (c4) [right=of c3] {};
    \node[circ] (c5) [right=of c4] {\Large $\substack{\!0, ~J-1\! \\ J+1}$};
    \node[circ] (c6) [right=of c5] {\Large $\substack{0,\ J \\ J+1}$};   
    \node[] (r2) [above=of c3] {};
    \node[crect] (r1) [above=1mm of r2] {$0, J+1$}; 
    \draw (c1) -- (r1);
    \draw (c2) -- (r1);
    \draw (c3) -- (r1);
    \draw (c5) -- (r1);
    \draw (c6) -- (r1);
    \node[rect] (s1) [below=3mm of c1] {$1$}; 
    \node[rect] (s2) [below=3mm of c2] {$2$};
    \node[rect] (s3) [below=3mm of c3] {$3$};    
    \node[rect] (s5) [below=3mm of c5] {$J-1$};    
    \node[rect] (s6) [below=3mm of c6] {$J$};      
    \draw[dotted] (c3) -- (c5);    
    \draw (c1) -- (s1); 
    \draw (c2) -- (s2);
    \draw (c3) -- (s3);
    \draw (c5) -- (s5);
    \draw (c6) -- (s6);
\end{tikzpicture}
}
\end{figure}

\begin{figure}[tb]
    \centering
\caption{
Modified junction tree for the graph of Wasserstein spline example from \cref{fig:wasserstein_spline}.
}
\label{fig:spline_modified_jt}
    \scalebox{0.7}{
\begin{tikzpicture}
	\tikzstyle{circ}=[circle, minimum size = 14mm, thick, draw =black!80, node distance = 3mm]
	\tikzstyle{rect}=[rectangle, minimum size = 6mm, thick, draw =black!80, fill=gray!50, node distance = 3mm]
 \tikzstyle{lrect}=[rectangle, minimum width=12mm, minimum height=8mm, thick, draw =black!80, node distance = 3mm]
  \tikzstyle{llrect}=[rectangle, minimum width=20mm, minimum height=8mm, thick, draw =black!80, node distance = 3mm]
\tikzstyle{inner}=[rectangle, draw=black!100, minimum size=5mm, fill=gray!50, inner sep=1pt]

	\node[circ] (c1) {$\substack{v_0,v_1 \\ x_0,x_1}$};
	\node[lrect] (12) [above right=of c1] {}; 
	\node[inner] at ([xshift=2mm] 12) {$x_1$}; 
        \node[left=-6mm of 12] {$v_1$}; 
 \node[circ] (c2) [below right=of 12] {$\substack{v_1,v_2 \\ x_1,x_2}$};
	\node[lrect] (23) [below right=of c2] {}; 
	\node[inner] at ([xshift=2mm] 23) {$x_2$}; 
    \node[left=-6mm of 23] {$v_2$}; 
	\node[circ] (c3) [above right=of 23] {$\substack{v_2,v_3 \\ x_2,x_3}$};
	\node[] (34) [right=of c3] {};
	\node[circ] (c4) [right=of 34] {$\substack{v_{J-2},v_{J-1} \\ x_{J-2},x_{J-1}}$};
	\node[llrect] (45) [below right=of c4] {}; 
	\node[inner] at ([xshift=4mm] 45) {$x_{J-1}$}; 
     \node[left=-10mm of 45] {$v_{J-1}$}; 
	\node[circ] (c5) [above right=of 45] {$\substack{v_{J-1},v_{J} \\ x_{J-1},x_J}$}; 
	\node[rect] (xN) [above right=of c5] {$x_{J}$};   
	\node[rect] (f2) [below left=of c1] {$x_0$}; 

	\draw (c1) -- (12);
	\draw (12) -- (c2);
	\draw (c2) -- (23);
	\draw (23) -- (c3);
	\draw[dotted] (c3) -- (c4);
	\draw (c4) -- (45);
	\draw (c5) -- (xN);
	\draw (c1) -- (f2); 
	\draw (12) -- (c2);
	\draw (23) -- (c3);
	\draw (45) -- (c5);
\end{tikzpicture}
}
\end{figure}

\subsection{Problem formulation}
We structure a coupled MOT problem \cref{eq:ot_graph_local}, which is equivalent to 
MOT problem \cref{eq:ot_graph},
on the junction tree in \cref{def:jt_general} by associating every cost term with a clique in $\cQ_C$.
 Moreover, we assume that each constraint is embedded in a separator in $\cQ_S$ as described in the following.
A constraint of the optimal transport problem can act on projections of the solution tensors $\bB_c$ onto a set of marginals. Recall that $\Gamma \subset \{1,2,\dots,J\} $ denote the set of index-sets corresponding to the constraints (see \cref{eq:omt_multi_discrete}). 
Then we assume that there exists a partition of $\Gamma$ such that for each $\bg$ in the partition there is a $c_\bg \in \cQ_S$ such that $\bg \subseteq c_\bg $, and 
any $c_\bg $ can contain only one constrained set $\bg$.
{Moreover, we denote the separators with constraints as $\cQ_\Gamma= \{c_\bg\} {\subseteq \cQ_S}$.}
Note that $\bg$ may not only contain one node, for example, if $ c_\bg = \{x_1,x_2, x_3\} $, we allow $\bg = \{x_1,x_3\}$.
Thus, we define the coupled MOT problem 
\begin{equation}
\begin{aligned} \label{eq:ot_graph_local}
 \min_{ \bB_c, c\in \cQ } & \sum_{c \in \cQ_C} \<\bC_{c} , \bB_{c} \> \\ \text{s.t. } \
 &  P_{{c_i} \cap {c_j}}( \bB_{c_i})  = P_{{c_i} \cap {c_j}}( \bB_{c_j}),\quad  \text{ for each } c_i \in \cQ_S, \text{ and all } c_j \in N(c_i), \\
       & P_\bg ( \bB_{c_\bg}) = \bmu_\bg, \quad  \text{for } 
       c_\bg \in \cQ_\Gamma {\subseteq \cQ_S},
   \qquad \< \bB_c, \mathbf{1} \> = 1, \quad  \text{for } c {\in \cQ_S \backslash \cQ_\Gamma.}
\end{aligned}
\end{equation}

The separators in $\cQ_S$ enforce that the cost cliques are consistent with each other, and thereby play a similar role as the separator sets in \cref{def:junction_tree}. 
To develop an algorithm for problem \cref{eq:ot_graph_local} we add a local regularization term for each clique that can contain a cost, i.e., for each $c\in \cQ_C$. Therefore, we generalize the expression for the reference joint distribution in \cref{eq:jointdist}. For 
$c\in \cQ_C$, let $c_j,c_k \in \cQ_S$ be the two neighbours, and define $ \bM_c = \bm_{c_j} \otimes \bm_{c_k}$, 
where  
$\bm_{c_i} = \bmu_\bg \otimes \mathbf{1} $ if $ \bg \subset c_i$ for some $\bg \subset \Gamma  $, and otherwise $\bm_{c_i} =\mathbf{1}$. 
We use the notation $\otimes \mathbf{1}$ to broadcast a tensor
to a bigger tensor of appropriate size. In the following, the broadcasting dimensions will always be obvious from the context.
For example, if ${c_i} = \{x_1,x_2,x_3\}, \bg = \{x_2\} ,$ then $\bm_{c_i} (x_1,x_2,x_3) = \bmu_\bg (x_2) $. 

Moreover, we add a mass constraint similar to \cref{eq:ot_graph_pairwise_mass}. 
The regularized problem thus reads
\begin{equation} \label{eq:ot_graph_local_reg}
		\begin{aligned}
			\min_{ \bB_c, c\in \cQ } & \sum_{c \in \cQ_C}  \< \bC_{c} , \bB_{c} \> + \epsilon \cH( \bB_c  | \bM_c  ) \\ \text{s.t. } \ 
			 &  P_{{c_i} \cap {c_j}}( \bB_{c_i}) = P_{{c_i} \cap {c_j}}( \bB_{c_j}),\quad  \text{ for each } c_i \in \cQ_S, \text{ and all } c_j \in N(c_i) \\
			                        & 
P_\bg ( \bB_{c_\bg }) = \bmu_\bg, \quad  \text{for } c_\bg \in \cQ_\Gamma \subseteq \cQ_S, \qquad
        \< \bB_c, \mathbf{1} \> = 1, \quad  \text{for } c \in \cQ_S \backslash \cQ_\Gamma.
		\end{aligned}
	\end{equation}

\subsection{Algorithm}
The problem in \cref{eq:ot_graph_local_reg}
 is similar to \cref{eq:ot_graph_pairwise_mass}, and we will extend the algorithm 
 developed in Section~\ref{sec:tree_local_reg} to handle general graphs.
More precisely,
following the derivations in \cref{sec:tree_local_reg}, we can derive a method for approximating the solution of \cref{eq:ot_graph_local_reg} similar to \cref{alg:tree_local_reg}.
Note that the second set of constraints in \cref{eq:ot_graph_local_reg} correspond to the constraints 
on $\Gamma$ in the tree-structured optimal transport problem \cref{eq:ot_graph_pairwise_mass}.
Special care needs to be applied to the first set of constraints.
In particular, in order to derive an efficient algorithm we need the neighbour sets in \cref{def:jt_general}.4 to be nested as described in the following assumption.
 
 \begin{assumption} \label{ass:inclusion}
 For all $ c \in \cQ_S 
 $ 
 with at least two neighbors, there is a inclusion sequence 
  $( c\cap c_{j_1} ) 
  =
  (c\cap c_{j_2} )\supseteq \dots \supseteq ( c \cap c_{j_\ell} ) $ 
that includes all $\ell$ neighbours $c_{j_1}, \dots ,c_{j_\ell} $ of $c$. 
In case a constraint set $\bg \subset \Gamma$ is embedded in the clique $c$, then this can be added after the smallest set in the inclusion, i.e., $( c\cap c_{j_\ell} ) \supseteq \bg$.
\end{assumption} 

\vspace{-0.6cm}
\begin{figure}[htb]
\noindent 
\begin{minipage}[c]{0.54\linewidth}
\begin{example}
We give an example of the nested structure in \Cref{ass:inclusion} in \cref{fig:nested}. Here, $c \!= \! \{x_1,x_3,x_4\}$, $c_1  \!= \! \{x_1,x_2,x_3,x_4\}$, $c_2  \!= \! \{x_1,x_3, x_4\}$, $c_3  \!= \! \{x_0,x_1, x_3\} $, and finally $\bg = \{x_3 \}.$
\end{example}
\end{minipage} \hfill
\begin{minipage}[c]{0.42\linewidth}
\centering
\captionof{figure}{Example on nested structure. }
\label{fig:nested}
\vspace{-10pt}
\scalebox{0.7}{
\begin{tikzpicture}
    \tikzset{
        circ/.style={circle, minimum size = 10mm, thick, draw=black!80, node distance=20mm},
        rect/.style={rectangle, draw=black!100, minimum width=14mm, minimum height=10mm},
        gray/.style={fill=gray!30},
        inner/.style={rectangle, draw=black!100, fill=gray!50, minimum size=5mm, inner sep=0.5pt}   
    }    
    \node[circ] (1) {$1,2,3, 4$};
    \node[rect, right=3mm of 1] (2) { };
    \node[inner] at ([xshift=0.5mm] 2.center) (i1) {$3$}; 
    \node[left=0.5mm of i1] {$1$}; 
    \node[right=0mm of i1] {$4$}; 
    \node[circ, right=3mm of 2] (3) {$1,3,4$};
    \node[circ, below=3mm of 2] (4) {$0, 1, 3$};
    \draw[thick] (1) -- (2);
    \draw[thick] (2) -- (3);
    \draw[thick, dashed] (3) -- ++(1,0);
    \draw[thick, dashed] (1) -- ++(-1,0);
    \draw[thick] (2) -- (4);
    \draw[thick, dashed] (4) -- ++(0,-1);
\end{tikzpicture}
}
\end{minipage}
\end{figure}
\begin{proposition} \label{prop:local_reg_dual}
The solution to \cref{eq:ot_graph_local_reg} is of the form 
$\bB_{c}  = \bK_{c} \odot \bU_{c} \odot 
\bM_{c}, $ for each $c \in \cQ_C$,
where $\bK_c = \exp(-\bC_c/\epsilon)$ and
$\bU_{c} = \left(\bu^{c_j}_{i_j}  \otimes \mathbf{1} \right) \odot \left( \bu^{c_k}_{i_k} \otimes \mathbf{1} \right) $ for two smaller tensors $\bu^{c_j}_{i_j}$ and $\bu^{c_k}_{i_k}$.
These tensors are dual variables associated with the two neighbors $c_j , c_k \in \cQ_S$ of $c$, and $i_j$, $i_k$ denote the position in the corresponding inclusion sequence as defined in Assumption~\ref{ass:inclusion}.
More precisely, denoting $\ell_c$ the length of the inclusion sequence for each $c\in \cQ_S$, the dual problem for \eqref{eq:ot_graph_local_reg} reads
\begin{equation}
\begin{aligned}
    \max_{
    \substack{\bu^c_i, i = 1,\ldots , \ell_c \text{ for } c \in \cQ_S \\ 
    \bu_{\bg}^c \text{ for } c_\bg \in \cQ_\Gamma \\
\rho_c \text{ for } c \in \cQ_S \backslash \cQ_\Gamma     
    }} & - \sum_{c \in \cQ_C}  \|\bK_{c} \odot \bU_c \odot \bM_c\|_1 
    + \! \sum_{c_\bg \in \cQ_\Gamma} \! \langle \bmu_\bg , \log ( \bu^c_{\bg }) \rangle 
    + \!\! \sum_{c \in \cQ_S \backslash \cQ_\Gamma} \!\! \rho_c
    \label{eq:whole_dual}  \\
    \mbox{subject to } &
 \sum_{i=1}^{\ell_c} \log(\bu^c_i \otimes \mathbf{1} )=\log(\bu^c_{\bg} \otimes \mathbf{1} ),  \quad  \text{for } 
 c \in \cQ_\Gamma 
     \\
    &   \sum_{i=1}^{\ell_c} \log(\bu^c_i \otimes \mathbf{1} )= \rho_c \mathbf{1} ,  \quad \text{for } c \in \cQ_S \backslash \cQ_\Gamma. 
\end{aligned}
\end{equation}
\end{proposition}
{The proof of \cref{prop:local_reg_dual} follows the standard Lagrangue dual theory.}
Recall that any separator set $c\in \cQ_S$ has a set of neighbours $c_1,\dots,c_\ell \in \cQ_C$ that define an inclusion sequence
$( c\cap c_{1} )   =  (c\cap c_{2} )\supseteq \dots \supseteq ( c \cap c_{\ell} )$ as in Assumption~\ref{ass:inclusion}.
To find the optimal dual variables of the dual problem \eqref{eq:whole_dual} we propose to iterate through all $c\in \cQ_S$, and update the dual variables $\bu^c_1,\dots \bu^c_{\ell_c}$.

Due to the fourth condition in \cref{def:jt_general}, the modified junction tree forms a bipartite graph, when we consider the separators in $ \cQ_S$ as nodes, and the cliques in $\cQ_C$ as edges. For example, in \cref{fig:wls_modified_jt}, the square node $\{ 0,J+1 \}$ alone is one set in a {\em bipartite} graph, and the bottom gray nodes belong to another set. 
The bipartite structure allows for a parallelized algorithm. In each update we can apply coordinate ascent to all the variables associated with separators in one set together.
Moreover, because those variables do not depend on each other, we can consider a subproblem of the dual, as we will describe in the following. Note that the iterative Bregman projection \cite{BenCarCut15} applies a similar idea to the Wasserstein barycenter problem.
To be precise, for the updates related to each $c\in \cQ_S$ it suffices to consider a subproblem of the dual \eqref{eq:whole_dual} , which for $c \in \cQ_\Gamma$ can be formulated as
\begin{equation}
    \max_{
    \substack{\bu_i, i = 1,\ldots , \ell \\ \bv_{\ell + 1}
    }} \sum_{i=1}^\ell - \< \bk_i , \bu_i \> 
    + \langle \bmu_\bg , \log ( \bv_{\ell + 1 }) \rangle 
  \quad  \mbox{s.t. } \sum_{i=1}^\ell \log(\bu_i \otimes \mathbf{1} )=\log(\bv_{\ell + 1} \otimes \mathbf{1} ). \label{eq:sub_dual}
\end{equation}

Here we have omitted the superscript of the dual variables. Moreover, we denote the dual variable $\bu_\gamma = \bv_{\ell+1}$, the reason for which will become apparent in the algorithm.
Also, we introduce the notation ${\bk_i := P_{c\cap c_i} (\bK_{c_i}  \odot \bM_{c_i} \odot (\bar\bu_i \otimes \mathbf{1}))}$, where $\bar\bu_i$ denotes the second dual variable defining the tensor $\bU_{c_i}$ together with $\bu_i$, as described in Proposition~\ref{prop:local_reg_dual}. Thus, we can write $\|\bK_{c_i} \odot \bU_{c_i} \odot \bM_{c_i}\|_1 = \< \bk_i , \bu_i \> $.
Note that the case $c\in \cQ_S \backslash \cQ_\Gamma$ is a special case of \eqref{eq:sub_dual}, where $\bv_{\ell+1} = {\exp(\rho_c)} \in \mathbb{R}$ is a scalar 
and $\bmu_\gamma=1$.

The optimal dual variables of the dual problem \eqref{eq:whole_dual} can now be found by iterating through all $c\in \cQ_S$, and solving a problem of the form \eqref{eq:sub_dual}. This can be done by
performing the steps described in the following.

\begin{proposition}\label{prop:general_algo}
The solution to \eqref{eq:sub_dual}, where $\ell>1$, can be found by performing the following steps.
\begin{enumerate}[noitemsep, leftmargin=*]
    \item Calculate geometric mean variables from the largest intersection cluster $c \cap c_1$ to the smallest $c \cap c_\ell$ as 
\begin{equation}
        \bq_2  = \left(\bk_1 \odot \bk_2 \right)^{1/2} \label{eq:q2_whole}; \qquad
    \bq_i  = \left( P_{c\cap c_i}( \bq_{i-1})^{i-1} \odot 
    \bk_i \right)^{1/i} , \qquad i=3,\dots, \ell  
\end{equation}
    \item Calculate intermediate variables from $\ba_\ell, \ldots, \ba_1$ and $\bv_\ell,\ldots, \bv_3$ from the smallest intersection cluster $c \cap c_\ell$ to the largest cluster $c \cap c_1$ as 
\begin{equation}
\begin{aligned}    
& \ba_{i }  = ( \bv_{i+1} {\otimes \mathbf{1}}) ^{1/i } \odot \bq_{i }, \qquad  \text{~~for } i = 2, \ldots, \ell \label{eq:av_updates} \\
& \bv_{\ell+1}  =  (\bmu_\gamma ./ P_\gamma(\bq_\ell))^\ell , \qquad \bv_{i } =(\ba_{i} ./ P_{c \cap c_{i}}(\bq_{i-1}))^{ i -1 },  \quad  \text{for } i = 3, \ldots, \ell 
\end{aligned}
\end{equation}
\item Update $\bu_1, \ldots, \bu_\ell$ as 
\begin{equation} \label{eq:u_updates_whole}
    \bu_{1} = \ba_2 ./ \bk_1 , \quad \bu_{i} = \ba_i ./ \bk_i  ~~\mbox{ for } 
    i={2,\dots,\ell}.
\end{equation}
\end{enumerate}
\end{proposition}

\Cref{prop:general_algo} only discusses the case of $\ell>1$. In the case $\ell=1$, the constraint in \eqref{eq:sub_dual} asserts that $\bu_1=\bv_2$. Thus, the objective in \eqref{eq:sub_dual} is optimized by $\bu_1= \bmu_\gamma ./ \bk_1$.

\begin{remark} \label{rem:more_neighbours}
    The method can be generalized to the case where the cost cliques have a larger number of neighbors, i.e., point 3. in \cref{def:jt_general} is generalized to: For all $c \in C_C$ it holds $|N(c)|\leq k$ for some $k \geq 2$.
    A parallelized algorithm for optimal transport problems on this type of junction tree can be set up by considering a $k$-partite partition, instead of the bipartite partition in the algorithm above.
\end{remark}

\begin{figure}[ht]
    \centering
\caption{Number of iterations for converging to the accuracy $\delta = 0.2$. The left column shows the number of iterations as a function of $d$
when $|E|$ is fixed ($|E| = 3$ in (a) and $|E| = 12$ in (b)),
and the right column vice versa with $d = 10$.}
    \label{fig:exp}
    \begin{subfigure}[b]{1\textwidth}
    \centering        
        \includegraphics[width=0.4\textwidth]{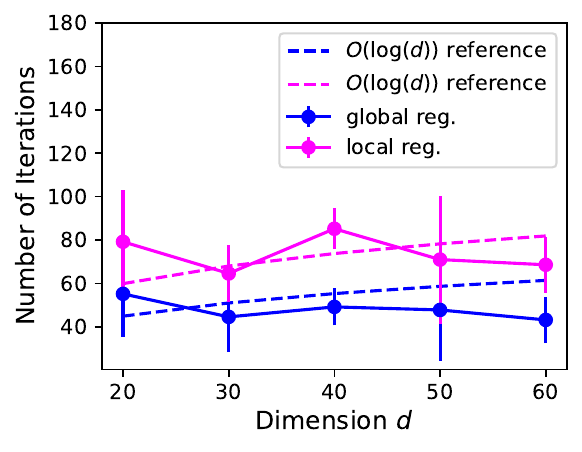}
        \includegraphics[width=0.4\textwidth]{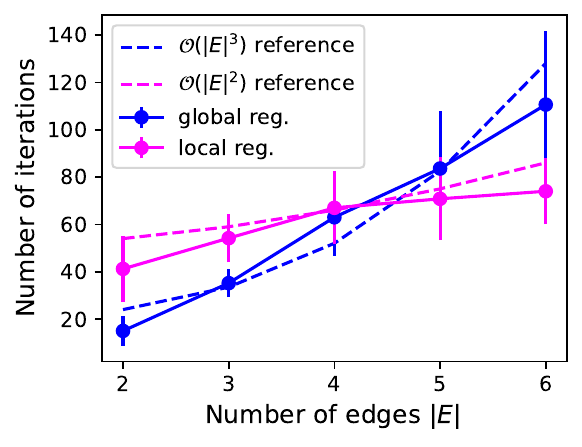}    
        \vspace{-15pt}
        \caption{Wasserstein Barycenter}
        \label{fig:barycenter}
    \end{subfigure}
    
    \begin{subfigure}[b]{1\textwidth}
    \centering      
        \includegraphics[width=0.4\textwidth]{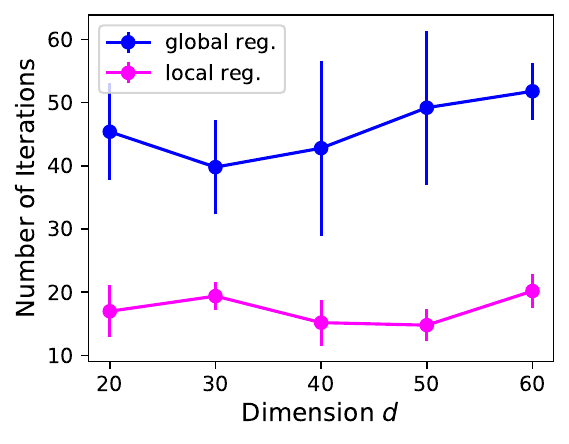}
        \includegraphics[width=0.4\textwidth]{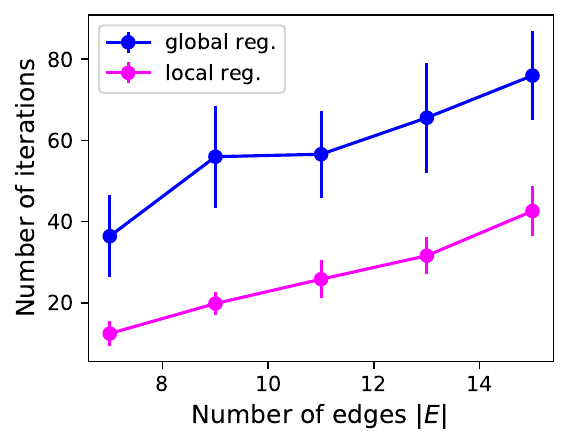}   
        \vspace{-15pt}
        \caption{Wasserstein least square}        
        \label{fig:wls}
    \end{subfigure}
\end{figure}

\section{Numerical examples}\label{sec:eg}

We present numerical results for two types of MOT problems. We examine the Wasserstein barycenter problem  as discussed in \cite[Example 1]{fan2021complexity} and structured according to the graph in \cite[Figure 1]{fan2021complexity}, and the Wasserstein Least Squares (WLS) example detailed in \cref{sec:wls}. The original graph and the modified junction tree for the WLS example are depicted in \cref{fig:wasserstein_least_square} and \cref{fig:wls_modified_jt}, respectively. For the sake of completeness, we also show modified junction tree for the Wasserstein spline example in \cref{fig:spline_modified_jt}. Notably, the barycenter problem is tree-structured, while the WLS problem forms a graph with a tree-width of two. \looseness=-1

For the barycenter problem, the cost matrices $C_{(j,k)}$ in \cref{eq:ot_graph_pairwise} are set as the squared Euclidean distances. In the WLS problem, the cost matrices follow \cref{eq:wls_cost} with coefficient $\alpha = 10$. The constrained marginal distributions $\{\mu_k\}_{k \in \Gamma}$, supported on a uniform grid of $d$ points between 0 and 1, are generated from the log-normal distribution and normalized to sum to one. We set the accuracy threshold $\delta = 0.2$ for the unregularized MOT problem. 

Our comparison involves \cref{alg:isbp_jt}, which solves a standard MOT problem with global regularization, and \cref{alg:tree_local_reg}, which addresses coupled bi-marginal OT with local regularization. We employ a random update rule for the global regularization method, consistent with the algorithm complexity analysis in \cite{fan2021complexity}. Our focus is on the number of iterations, anticipating that the parallel mechanism introduced by the local regularization method could reduce them. Each experiment is repeated five times with different random seeds, and the average number of iterations is reported in \cref{fig:exp}. Theoretical complexity bounds for the barycenter problem are also illustrated as dashed lines, with bounds for global and local regularization derived from \cite[Theorem 2]{fan2021complexity} and \cref{thm:tree_local_reg_complexity}, respectively. 

From the observations in \cref{fig:barycenter}, it is evident that the experimental results are consistent with the theoretical bounds. Notably, local regularization demonstrates the potential to reduce the number of iterations by a factor of $|E|$, leading to a lower iteration count beyond a certain threshold $|E|$. Furthermore, both global and local regularization exhibit a logarithmic dependency on the dimension, with a $\log(d)$ relationship. In contrast, for the results presented in \cref{fig:wls}, we have not established theoretical bounds for local regularization within this paper; this constitutes an area for future investigation.

\section{Conclusion}

In this paper, we study the MOT problem with graph-structured costs. We propose a novel approach to solve the MOT problem using local entropy regularization. Specifically, we reformulate the entire MOT problem as a sum of smaller MOT problems. When the underlying graph is a tree, these smaller problems involve dual variables that can be divided into two groups, forming a bipartite graph. We develop an efficient parallel algorithm that alternately updates these two groups of dual variables. A key advantage of this approach is its parallelizability, making it well-suited for GPU implementation. We also provide a rigorous complexity analysis of this approach and extend it to handle MOT problems associated with general graphs.

\section{Proofs of Theorems}



\begin{proof}[Proof of \cref{thm:tree_dual_var}]
    We perform a block-coordinate descent in \cref{eq:tree_dual_mass_consv}. First, minimizing \cref{eq:tree_dual_mass_consv} with respect to $u_{(j,k_j)}$ for $j\in \Gamma$, while keeping the other variables fixed, yields the first set of updates.
    To minimize the objective with respect to $\rho_j$ and $u_{(j,k)}$ for $j \notin \Gamma$, relax the constraint with Lagrangian multiplier $\beta$, which gives the Lagrangian $L(\{\rho_j\}, \{u_{(j,k)}\}, \{\beta_j\} =$
    \begin{equation} \label{eq:tree_lagrangian_alpha}
         \!\! \sum_{(j,k)\in E} \!\! \|B_{(j,k)}\|_1 - \sum_{j \in \Gamma} \log (u_{(j,k_j)})^\top \mu_{j}   - \!\! \sum_{j \in V \backslash \Gamma} \!\! \left( \rho_j + \beta_j^\top \! \left( \sum_{k \in N(j)} \!\! \log(u_{(j,k)}) - \rho_j \mathbf{1} \right) \!\! \right). \!\!
    \end{equation}
    The Lagrangian is minimized with respect to $u_{(j,k)}$ for $j \notin \Gamma$ if
    \begin{equation} \label{eq:tree_v_alpha}
    \beta_j = u_{(j,k)} \odot (K_{(j,k)} (u_{(k,j)} \odot \gamma_k )),\quad \text{ for } k \in N(j).    
    \end{equation}
    Multiplying for all $k \in N(j) $ yields
    $\bigodot_{k \in N(j)} \beta_j = e^{\rho_j} \bigodot_{k \in N(j)} \left( K_{(j,k) } (u_{(k,j)} \odot \gamma_k ) \right)$. 
   Denote 
   $v_j = \left( \bigodot_{k \in N(j) }  (K_{(j,k)} (u_{(k,j)} \odot \gamma_k )) \right)^{1/|N(j)|}, $
   it follows that 
    $\beta_j = e^{\rho_j/|N(j)| } v_j,$
    which together with \cref{eq:tree_v_alpha} leads to the updates for $u_{(j,k)}$, $j \notin \Gamma$.
    Note also that the Lagrangian \cref{eq:tree_lagrangian_alpha} is only bounded in $\rho_j$ if it holds $\beta_j^\top \mathbf{1} = 1$.
    Thus, $e^{\rho_j/ |N(j)|}  v_j^\top \ett = 1$ and the update for $\rho_j$ follows. 
\end{proof}

\begin{proof}[Proof of \cref{thm:num_iter}]
Define $\Delta_t:=     f(u^t,\rho^t) - f( u^*, \rho^*). $ Then according to the \cref{lem:f_upper_bound_tree}, \cref{lem:f_lower_bound_tree}, there is
$\Delta_t - \Delta_{t+1} \ge \frac{1}{22|E|} \left( \max\left\{ \frac{\epsilon \Delta_t}{2 C_\infty}, \delta' \right\} \right)^2. $
On the one hand, we have
$\frac{\epsilon^2 \Delta_{t+1}}{88|E| C_\infty^2} \le \frac{\epsilon^2 \Delta_{t}}{88 |E|  C_\infty^2} - \left(\frac{\epsilon^2 \Delta_t}{88 |E|  C_\infty^2} \right)^2. $
Note that $0 \le a_{t+1} \le a_t -a_t^2$ implies $0 \le a_t \le 1$. This implies whenever $\Delta_{t+1} \ge 0$, we have $0 \le \frac{\epsilon^2 \Delta_{t}}{88|E|C_\infty^2}  \le 1$ for any $t >1$.
Also note that $0 \le a_{t+1} \le a_t -a_t^2$ implies ${1}/{a_{t+1}} - {1}/{a_{t}} \ge 1$ and further implies ${1}/{a_t} - {1}/{a_1} \ge t-1$.
Thus,
\begin{align}\label{eq:time_bound1_tree}
t \le  1+ \frac{88|E|C_\infty^2}{\epsilon^2 \Delta_{t}} - \frac{88|E|C_\infty^2}{\epsilon^2 \Delta_{1}}.
\end{align}
On the other hand, 
\begin{align}\label{eq:time_bound2_tree}
\Delta_{t+\tau} \le \Delta_{t } - \frac{\delta'^2 \tau}{22 |E|}, \quad t, \tau \ge 0.
\end{align}
Fixing an number $\phi>0$, and considering using \cref{eq:time_bound1_tree} for $\Delta_t \ge \phi$ and  using \cref{eq:time_bound2_tree} for $\Delta_t < \phi$. We minimize the sum of two estimates over the parameter $\phi \in (0, \Delta_1]$.
Combining the above two inequalities, we get the total number of iterations satisfies
\begin{align}
    t \le \min_{0 \le \phi \le \Delta_1}
    \left( 1 + \frac{88|E|C_\infty^2}{\epsilon^2 \phi} - \frac{88|E|C_\infty^2}{\epsilon^2 \Delta_1} + 1+ \frac{22 |E| \phi }{\delta'^2} \right) \le 2+ \frac{88 |E| C_\infty}{\epsilon \delta'}.
\end{align}
\end{proof}

\begin{proof}[Proof of \cref{thm:tree_local_reg_complexity}]
With the specific choices $\epsilon = \frac{\delta}{ 4 |E| \log(d) } $ and $\delta' = \frac{\delta}{ 8 C_\infty} $,  \cref{lem:round_applied} gives
\begin{align}
 \sum_{(j,k) \in E } \< C_{(j,k)}, \widehat B_{(j,k)} \> - \sum_{(j,k) \in E } \< C_{(j,k)},  B^*_{(j,k)} \> 
\le \delta.
\end{align}
By \cref{thm:num_iter}, the stopping time $t$ satisfies $t \le 2+ \frac{88|E|C_\infty}{\delta' \epsilon} = \mathcal O \left( \frac{|E|^2  C_\infty^2 \log(d)}{\delta^2 } \right). $ Since each iteration of \cref{alg:tree_local_reg} requires $\mathcal O(|E|d^2)$ operations, a solution of 
\cref{alg:tree_local_reg} is achieved in 
$\mathcal O \left( \frac{|E|^3  C_\infty^2 d^2 \log(d)}{\delta^2 } \right) . $
 On the other hand, \cref{alg:round} takes $\mathcal O(|E|d^2 )$ (see Lemma 7 in \cite{altschuler2017near}). Therefore, the bound on total complexity
follows.
\end{proof}

\section{Proofs of Propositions}


\begin{proof}[Proof of \Cref{cor:equi}]
Firstly, when $j \in \Gamma$, we have $\gamma_j = \mu_j $, and
\begin{equation*}
u^{t+1}_{(j,k_j)} \overset{\cref{eq:fixed_node_update_mass_cons}}{ = }  \mathbf{1} ./ (K_{(j,k_j)} u^t_{(k_j,j)}) 
 \overset{ \cref{eq:bi-marginal_B_sol} }{=} \mathbf{1} ./ ( (B^t_{(j k_j)} \ett) ./ (u^t_{(j, k_j)} \odot   \mu_j ) )  \overset{ \cref{eq:q} }{=} u^t_{(j, k_j)} \odot   \mu_j ./  q^t_{(j,k_j)} .
\end{equation*}
When $j \notin \Gamma$, we have $\gamma_j = \ett$. It follows that
\begin{align*}
v^{t+1}_j & \overset{ \cref{eq:free_node_update_mass_cons} }{=} \left( \bigodot_{k\in N(j)} \left( K_{(j,k) } (u^t_{(k,j)} 
   \odot \gamma_k ) \right)
 \right)^{\frac{1}{|N(j)|}}   \overset{ \cref{eq:bi-marginal_B_sol} }{=} \left( \bigodot_{k\in N(j)} \left( B^t_{(j,k)} \ett ./ u^t_{(j,k)} 
    \right)
 \right)^{\frac{1}{|N(j)|}}  \\
 & = \left( \bigodot_{k\in N(j)} q^t_{(j,k)}   ) \right)^{{\frac{1}{|N(j)|}} } / e^{\rho_j^t /|N(j)| } = q^{t+1}_{j} / e^{\rho_j^t /|N(j)| } ,
\end{align*}
where we use the constraint $\sum_{k \in N(j)} \log(u_{(j,k)}) = \rho_j \mathbf{1}$ of the dual problem \cref{eq:tree_dual_mass_consv} in the third equality, and the definition of $q_j^{t+1}$ in the last equality. 
This implies that
\begin{align*}\label{eq:rho_update}
\rho^{t+1}_j 
 = - |N(j)| \log \left(\mathbf{1}^\top v^{t+1}_j \right) 
 = \rho^t_j - |N(j)| \log(\ett^\top q^{t+1}_{j}) =   \rho^t_j - |N(j)| \log(\|q^{t+1}_{j}\|_1 ).
\end{align*}
Furthermore,
\begin{align*}
  u^{t+1}_{(j,k)}
\overset{ \cref{eq:free_node_update_mass_cons} }{=} & e^{\rho^{t+1}_j /|N(j)|}  v^{t+1}_j ./ (K_{(j,k)} (u^t_{(k,j)} \odot \gamma_k))   \\
\overset{ \cref{eq:bi-marginal_B_sol} }{=}   
& e^{(\rho^{t+1}_j -\rho^{t}_j ) /|N(j)|} q^{t+1}_{j}  ./ (B^t_{(j,k)} \ett ./ u^t_{(j,k)})  
= u_{(j,k)}^t 
\odot \frac{q_j^{t+1}}{\|q_j^{t+1}\|_1}  ./  q_{(j,k)}^t .
\end{align*}
\end{proof}

\begin{proof}[Proof of \Cref{prop:general_algo}]
We introduce auxiliary variables ${\bv_i}$, $i=1,\dots,\ell$, and rewrite the constraints \cref{eq:sub_dual}, resulting in the following reformulation of the problem 
\begin{align}
    \max_{
    \substack{\bu_i, i = 1,\ldots , \ell \\ \bv_i, i = 3,\ldots , \ell + 1
    }} & - \sum_{i=1}^\ell \< \bk_i , \bu_i \> 
    + \langle \bmu_\bg , \log (\bv_{\ell + 1 }) \rangle  \label{eq:dual1} \\
    \mbox{subject to } 
    & \log(\bu_1 ) + \log(\bu_2 ) = \log(\bv_3 \otimes \mathbf{1} ), \label{eq:uv1} \\
    & \log(\bv_i ) + \log(\bu_i  ) = \log(\bv_{i+1} \otimes \mathbf{1} ), \quad \mbox{ for } i=3, \ldots, \ell -1 \\
    & \log(\bv_\ell ) + \log(\bu_\ell ) = \log(\bv_{\ell +1} \otimes \mathbf{1} ) . \label{eq:uv_last}
\end{align}
Note here that $\bv_i $ and $\bu_i$ have the same dimension.
The Lagrange relaxed problem of \cref{eq:dual1} is 
\begin{equation}
\begin{aligned}
    & \max_{\substack{
     \bu_i, i = 1,\ldots , \ell  \\ \bv_i, i = 3,\ldots , \ell+ 1 
    }}  - \sum_{i=1}^\ell  \< \bk_i , \bu_i \>  + \langle \bmu_\bg , \log (\bu_{\ell + 1 }) \rangle  
     +\langle \ba_2 , \log(\bu_1)+\log(\bu_2 )-\log(\bv_3 \! \otimes \! \mathbf{1} )\rangle\\
    & \ + \!\sum_{i=3}^{\ell -1}\langle \ba_i, \log(\bv_i)+\log( \bu_i  )-\log(\bv_{i+1} \!\otimes \! \mathbf{1} )\rangle +\langle \ba_\ell , \log(\bv_\ell )+\log(\bu_\ell )-\log(\bv_{\ell + 1} \! \otimes \! \mathbf{1} )\rangle,
\end{aligned}
\end{equation}
where $\ba_i$ are newly introduced 
Lagrange variables.
Maximizing with respect to $\bv_i$ we get $ \ba_i = P_{c\cap c_i}(\ba_{i-1})$ for $i=3,\ldots, \ell$. 
Maximizing with respect to $\bv_{\ell + 1}$, we get that $
P_\bg( \ba_\ell) =\bmu_\bg.$ Further, maximizing with respect to $\bu_i$, for $i=1,\dots,\ell$, we get \cref{eq:u_updates_whole}.
To satisfy the constraint \eqref{eq:uv1} we need $\bu_1 \odot \bu_2 = \bv_3 \otimes \mathbf{1}$. Plugging in the expressions for $\bu_1$ and $\bu_2$ in \eqref{eq:u_updates_whole}, we get the expression for $\ba_2$ in \eqref{eq:av_updates} together with \eqref{eq:q2_whole}.
The update of $\bv_3$ in \eqref{eq:av_updates} follows from $\ba_3 = P_{c \cap c_i}(\ba_2)$.
Similarly, the remaining updates for $\bq_i$, $\ba_i$, and $\bv_i$, for $i=3,\dots,\ell$, follow from the constraints $\bv_i \odot \bu_i = \bu_{\ell+1} \otimes \mathbf{1}$ and $\ba_i = P_{c\cap c_i}(\ba_{i-1})$.
Finally, the expression for $\bv_{\ell+1}$ is a consequence of $\bmu_\gamma = P_{\gamma}(\ba_\ell)$.
\end{proof}

\section{Proofs of Lemmas}

\begin{proof}[Proof of \Cref{cor:tree_mass}]
When $j \in \Gamma \bigcap S$, the updates of $u_{(j,k)}$ keep 
$\|B^{t+1}_{(j,k_j)}\|_1 =1$, since it matches the projection $B^{t+1}_{(j,k_j)} \mathbf{1}$ and the true marginal $\mu_j$. 
When $j \in (V \backslash \Gamma) \bigcap S$, by \cref{cor:equi}, 
the updates of $u_{(j,k)}$ are equivalent to \cref{eq:free_node_update_mass_cons} in \cref{thm:tree_dual_var}, thus they satisfy \cref{eq:tree_v_alpha}, which lead to 
$\|B_{(j,k)}^{t+1} \| = \mathbf{1}^\top \beta^{t+1}_j. $ Since $\mathbf{1}^\top \beta^{t+1}_j=1$ holds, we get $\|B_{(j,k)}^{t+1} \|=1.$
\end{proof}

\begin{proof}[Proof of \Cref{lem:tree_dual_varia_bounds}]
By \cref{eq:fixed_node_update_mass_cons}, for $j \in \Gamma$, following the logic in \cite[Lemma~1]{dvurechensky2018computational}, 
there is
\begin{align}
 \max_{i=1,\dots,d} (\lambda_{(j,k)}^t)_i - \min_{i=1,\dots,d} (\lambda_{(j,k)}^t)_i  \leq 
 \|C_{(j,k)}\|_\infty.
\end{align}
By \cref{eq:free_node_update_mass_cons}, for $j \in V \backslash  \Gamma $, following \cite[Lemma~4]{kroshnin2019complexity}, there is
\begin{align}
 \max_{i=1,\dots,d} (\lambda_{(j,k)}^t)_i - \min_{i=1,\dots,d} (\lambda_{(j,k)}^t)_i  
 \leq 
 \|C_{(j,k)}\|_\infty + \frac{1}{|N(j)|} \sum_{i \in N(j)} \|C_{(j,i)}\|_\infty.
\end{align}
Combining these results gives the bound.
\end{proof}

\begin{proof}[Proof of \Cref{lem:f_upper_bound_tree}]
 By the gradient inequality for convex functions it holds that
 \begin{align*} 
      &  f(u^t,\rho^t) - f(u^*,\rho^*) = f(\lambda^t,\rho^t) - f(\lambda^*,\rho^*) \\
      & \leq  
            \langle \nabla f(\lambda^t,\rho^t) , \left( (\lambda^t,\rho^t) - (\lambda^*,\rho^*) \right)^\top \rangle \\
            & = \! \sum_{j \in \Gamma}  \langle \lambda_{(j,k_j)}^t \! - \! \lambda_{(j,k_j)}^*, \frac{1}{\epsilon} B_{(j,k_j)}^t \mathbf{1} \! - \! \frac{1}{\epsilon} \mu_j \rangle \!   + \!\sum_{j \notin \Gamma} \sum_{k \in N(j)} \!\!\! \langle \lambda_{(j,k)}^t \! - \! \lambda_{(j,k)}^* , \frac{1}{\epsilon} q_{(j,k)}^t \rangle  - \! \sum_{j \notin  \Gamma} ( \rho_j^t \! - \! \rho_j^*) \\
            & = \frac{1}{\epsilon} \sum_{j \in \Gamma}  \langle \lambda_{(j,k_j)}^t - \lambda_{(j,k_j)}^*,  B_{(j,k_j)}^t \mathbf{1} - \mu_j \rangle  \\
            & \quad + \sum_{j \notin \Gamma} \sum_{k \in N(j)} \left\langle \frac{1}{\epsilon} \lambda_{(j,k)}^t - \frac{1}{|N(j)|} \rho_j^t \mathbf{1} - \frac{1}{\epsilon} \lambda_{(j,k)}^* 
           - \frac{1}{|N(j)|}  \rho_j^* \mathbf{1}  ~,~  q_{(j,k)}^t \right\rangle  .     \label{eq:proof_fdistance_tree}
\end{align*}     
    Here, we have used that $\langle \mathbf{1}, q_{(j,k)}^t \rangle = 1$, and thus $ - \rho_j^t + \rho_j^* = \langle - \rho_j^t \mathbf{1} + \rho_j^* \mathbf{1}, q_{(j,k)}^t \rangle $. 
 We bound the first sum as in \cite[proof of Lemma~2]{fan2021complexity} (compare also \cite[proof of Lemma~2]{dvurechensky2018computational}),
 which leads to
     \begin{equation} \label{eq:proof_fixed_dual}
             \sum_{j \in \Gamma}  \langle \lambda_{(j,k_j)}^t - \lambda_{(j,k_j)}^*, 
             B_{(j,k_j)}^t \mathbf{1} - 
             \mu_j \rangle 
             \leq 2 C_\infty \sum_{j \in \Gamma}  \|B_{(j,k_j)}^t \mathbf{1} - \mu_j\|_1.
     \end{equation}
     To bound the second sum in the final expression of \cref{eq:proof_fdistance_tree} note that for any $j \notin \Gamma,$
 $$ \sum_{k \in N(j)} \left\langle \frac{1}{\epsilon} \lambda_{(j,k)}^t - \frac{1}{|N(j)|} \rho_j^t \mathbf{1} 
           ~,~  \bar q_j^t \right\rangle 
 =  \left\langle \sum_{k \in N(j)} \frac{1}{\epsilon} \lambda_{(j,k)}^t -  \rho_j^t \mathbf{1} 
           ~,~  \bar q_j^t \right\rangle  = 0 , $$
 since $\sum_{k\in N(j)} \lambda_{(j,k)}^t/\epsilon =\rho_j^t \mathbf{1}$. The same equality holds with $\lambda_{(j,k)}^*$ and $\rho_j^*$.
     Moreover, note that it holds
    $$ \sum_{k\in N(j)} \langle  \mathbf{1}, q_{(j,k)}^t - \bar q_j^t \rangle = \left\langle \mathbf{1} , \sum_{k\in N(j)} q_{(j,k)}^t - |N(j)| \bar q_j^t \right\rangle = 0 .$$
     Using these relations it follows
      \begin{align} 
& ~~~\sum_{k \in N(j)} \left\langle \frac{1}{\epsilon} \lambda_{(j,k)}^t - \frac{1}{|N(j)|} \rho_j^t \mathbf{1} - \frac{1}{\epsilon} \lambda_{(j,k)}^* 
           - \frac{1}{|N(j)|}  \rho_j^* \mathbf{1}  ~,~  q_{(j,k)}^t \right\rangle \\
            &= \frac{1}{\epsilon} \sum_{k \in N(j)} \left\langle  \lambda_{(j,k)}^t-  \lambda_{(j,k)}^* 
             ~,~  q_{(j,k)}^t -\bar q_j^t \right\rangle  
             \leq \frac{2 C_\infty}{\epsilon} \sum_{k \in N(j)} \|B_{(j,k)}^t \mathbf{1} - \bar q_j^t\|_1. \label{eq:proof_free_dual}
    \end{align}
Here, the last step follows with the same approach as for \cref{eq:proof_fixed_dual} (see also \cite[proof of Lemma~5]{kroshnin2019complexity}).
Finally, combining \cref{eq:proof_fdistance_tree}, \cref{eq:proof_fixed_dual} and \cref{eq:proof_free_dual} completes the proof.
\end{proof}

\begin{proof}[Proof of \Cref{lem:f_lower_bound_tree}]
Using the relations in \cref{alg:tree_local_reg} and \cref{cor:tree_mass}, we get
    \begin{align*}
    & f(u^t,\rho^t) - f(u^{t+1}, \rho^{t+1}) 
    = \sum_{j \in \Gamma} (\log u_{(j,k_j)}^{t+1} - \log u_{(j,k_j)}^t )^\top \mu_j 
    + \sum_{j \notin \Gamma} (\rho_j^{t+1} - \rho_j^t ) \\
    =& \sum_{j \in \Gamma} \left( {\cH} \left( \mu_j | B_{(j,k_j)}^t \mathbf{1} \right) + 1\right)  - \sum_{j \notin \Gamma} |N(j)| \log \left(\|q_j^{t+1}\|_1 \right).
    \end{align*}
    By Pinsker's inequality the first sum can be bounded as
    $$  \sum_{j \in \Gamma} \left( {\cH} \left( \mu_j | B_{(j,k_j)}^t \mathbf{1} \right) + 1\right) 
    \geq \frac{1}{2}  \sum_{j \in \Gamma} \|B_{(j,k_j)}^t \mathbf{1} - \mu_j\|_1^2 .$$
    Moreover, since $\|q_j^{t+1}\|_1 \leq 1$ and $\|q_{(j,k)}^t\|_1 = 1$, we can bound
        \begin{align*}
    & - \log \left(\|q_j^{t+1}\|_1 \right) \geq  1 - \|q_j^{t+1}\|_1 = \frac{1}{|N(j)|}\sum_{k \in N(j)} \| q_{(j,k)}^t \|_1 - \|q_j^{t+1}\|_1 \\
    = & \mathbf{1}^\top \left( \frac{1}{|N(j)|} \sum_{k \in N(j)} q_{(j,k)}^t  - q_j^{t+1} \right) \geq \frac{1}{11 |N(j)|} \sum_{k \in N(j)} \|q_{(j,k)}^t - \bar q_j^t\|_1^2, 
    \end{align*}
    where the last inequality is proved in \cite[proof of Lemma~6]{kroshnin2019complexity}.
Therefore, 
    \begin{align*}
      f(u^t,\rho^t) - f(u^{t+1}, \rho^{t+1}) 
  & \ge  \frac{1}{2}  \sum_{j \in \Gamma} \|B_{(j,k_j)}^t \mathbf{1} - \mu_j\|_1^2 
    + \frac{1}{11} \sum_{j \notin \Gamma}  \sum_{k \in N(j)} \|q_{(j,k)}^t - \bar q_j^t\|_1^2 \\
       & \ge  \frac{1}{11} \left(\sum_{j \in \Gamma} \|B_{(j,k_j)}^t \mathbf{1} - \mu_j\|_1^2  +  \sum_{j \notin \Gamma}  \sum_{k \in N(j)} \|q_{(j,k)}^t - \bar q_j^t\|_1^2  \right) \\
    & \ge  \frac{1}{22|E|} \left(\sum_{j \in \Gamma} \|B_{(j,k_j)}^t \mathbf{1} - \mu_j\|_1  +  \sum_{j \notin \Gamma}  \sum_{k \in N(j)} \|q_{(j,k)}^t - \bar q_j^t\|_1   \right)^{\!2} \! ,     
    \end{align*}
    where we use Cauchy-Schwarz in the final inequality.
\end{proof}

\begin{proof}[Proof of \Cref{lem:round}]
Denote the partition where the nodes have mismatch as $S$. Without loss of generality, we assume $j \in S$ and $k \in V \backslash S$. 
 \begin{align*}
& \sum_{(j,k) \in E } \< C_{(j,k)}, B_{(j,k)} \> - \sum_{(j,k) \in E } \< C_{(j,k)}, \widehat B_{(j,k)} \>
 = \sum_{(j,k) \in E } \< C_{(j,k)}, B_{(j,k)} - \widehat B_{(j,k)} \>  \\
& = \sum_{j \in \Gamma } \< C_{(j,k_j)}, B_{(j,k_j)} - \widehat B_{(j,k_j)} \> + \sum_{j \notin \Gamma } \sum_{k \in N(j) } \< C_{(j,k)}, B_{(j,k)} - \widehat B_{(j,k)} \> .
 \end{align*}
By \holder inequality and Lemma 7 in \cite{altschuler2017near}, the first part can be bounded and 
\begin{align*}
\sum_{j \in \Gamma } \< C_{(j,k_j)}, B_{(j,k_j)} - \widehat B_{(j,k_j)} \>
  &\le   \sum_{j \in \Gamma }   \| C_{(j,k_j)}\|_\infty \| B_{(j,k_j)}  - \widehat B_{(j,k_j)} \|_1   \\
 &\le 2 C_\infty \sum_{j \in \Gamma }    \| \mu_j - B_{(j,k)} \ett \|_1 .
\end{align*}
Similarly, for the second part,
\begin{align*}
& \sum_{j \notin \Gamma } \sum_{k \in N(j) } \< C_{(j,k)}, B_{(j,k)} - \widehat B_{(j,k)} \> 
  \le   \sum_{j \notin \Gamma }  \sum_{k \in N(j) } \| C_{(j,k)}\|_\infty \| B_{(j,k)}  - \widehat B_{(j,k)} \|_1 \\
 & \le  
   2 \sum_{j \notin \Gamma } \sum_{k \in N(j) }   \| C_{(j,k)}\|_\infty \| \bar q_j - B_{(j,k)} \ett \|_1 
   \le  
   2 C_\infty \sum_{j \notin \Gamma } \sum_{k \in N(j) }    \| \bar q_j - B_{(j,k)} \ett \|_1 .   
\end{align*}
Summing them up gives
\begin{align*}
 & \sum_{(j,k) \in E } \< C_{(j,k)}, B_{(j,k)} \> 
- \sum_{(j,k) \in E } \< C_{(j,k)}, \widehat B_{(j,k)} \>  \\
\le & 2 C_\infty  \left( \sum_{j \in \Gamma }   \| \mu_j - B_{(j,k_j)} \ett) \|_1  + \sum_{j \notin \Gamma } \sum_{k \in N(j) }   \| \bar q_j (B) - B_{(j,k)} \ett \|_1 \right) .
\end{align*}
It follows the exact same way to obtain the same bound for $ \sum_{(j,k) \in E } \< C_{(j,k)}, \widehat B_{(j,k)} \> -  \sum_{(j,k) \in E } \< C_{(j,k)}, B_{(j,k)} \> $.
\end{proof}

\begin{proof}[Proof of \Cref{lem:round_applied}]
Let $\{\widehat{Y}_{(j,k)}\}$ denote the output from \cref{alg:round} with inputs $\{ B^*_{(j,k)}\}$ and $\{\widetilde B_{(j,k_j)} \ett \}_{j \in \Gamma} $. Note that $\widetilde B_{(j,k)}$ is the optimal solution to 
\begin{align*}
\min_{B \in \Pi(\widetilde B_{(j,k)} \ett, ~ \widetilde B^\top_{(j,k)} \ett ) } \< C_{(j,k)} , B \>  + \epsilon \cH(B | M_{(j,k)} ), \quad \text{for } \forall (j,k) \in E
\end{align*}
which can be verified by checking KKT conditions. Thus it holds that 
\begin{align*}
\< C_{(j,k)}, \widetilde B_{(j,k)} \> + \epsilon \cH(\widetilde B_{(j,k)} | M_{(j,k)} ) \le \< C_{(j,k)}, \widehat Y_{(j,k)} \>   + \epsilon \cH(\widehat Y_{(j,k)} | M_{(j,k)} ).
\end{align*}
Since $\< \widetilde B_{(j,k)}, \log(\widetilde B_{(j,k)}) \> \ge - 2 \log(d)$ and $\< \widehat Y_{(j,k)}, \log(\widehat Y_{(j,k)}) \> \le 0 $ it follows that 
\begin{align*}
& \< C_{(j,k)}, \widetilde B_{(j,k)} \> -   \< C_{(j,k)}, \widehat Y_{(j,k)} \> 
 \le \epsilon \cH(\widehat Y_{(j,k)} | M_{(j,k)} ) - \epsilon \cH(\widetilde B_{(j,k)} | M_{(j,k)} )   \\ 
& \le - \epsilon \< \widetilde B_{(j,k)}, \log(\widetilde B_{(j,k)}) \> +  \epsilon \< \widetilde B_{(j,k)} - \widehat Y_{(j,k)}  , \log M_{(j,k)}  \>  \\
& \le 2\epsilon \log(d) + \epsilon \<\log \gamma_j , \widetilde B_{(j,k)} \ett - \widehat Y_{(j,k)} \ett \> + \epsilon \<\log \gamma_k , \widetilde B_{(j,k)}^\top \ett - \widehat Y_{(j,k)}^\top \ett\>  
= 2\epsilon \log(d) , \label{eq:CB1}
\end{align*}
where the last equality holds because $\widetilde B_{(j,k)} \ett = \widehat Y_{(j,k)} \ett $ whenever $j\in \Gamma$, and $\gamma_j = \ett$ whenever $j \notin \Gamma$.

The plans $\{ B^*_{(j,k)}\}$ already satisfy the condition $\bar q_j(B^*) = B^*_{(j,k)} \ett $ for $\forall j \in \Gamma, k \in N(j)$.
Thus \cref{lem:round} gives 
\begin{equation}
\sum_{(j,k) \in E } \< C_{(j,k)}, \widehat Y_{(j,k)} \> 
- \sum_{(j,k) \in E } \< C_{(j,k)},  B^*_{(j,k)} \>  \le  2 C_\infty  \left( \sum_{j \in \Gamma }   \| \mu_j - \widetilde B_{(j,k_j)} \ett) \|_1  
\right)  , \label{eq:CB2}
\end{equation}
\begin{equation}
\begin{aligned}
& \sum_{(j,k) \in E } \< C_{(j,k)}, \widehat B_{(j,k)} \> 
- \sum_{(j,k) \in E } \< C_{(j,k)},  \widetilde B_{(j,k)} \> \\
& \qquad \le  2 C_\infty  \left( \sum_{j \in \Gamma }   \| \mu_j - \widetilde B_{(j,k_j)} \ett) \|_1   + \sum_{j \notin \Gamma } \sum_{k \in N(j) }   \| \bar q_j (\widetilde B) - \widetilde B_{(j,k)} \ett \|_1 \right) . \label{eq:CB3}
\end{aligned}
\end{equation}
Summing up \cref{eq:CB1}, \cref{eq:CB2} and \cref{eq:CB3} concludes the proof.
\end{proof}


\bibliographystyle{siamplain}
\bibliography{references}

\newpage

\begin{center}
{ \Huge \bf Supplementary material }
\end{center}

\vspace{20pt}

In this supplementary material, we present a concise framework for MOT problems with global entropy regularization. 
Extending our previous work \cite{fan2021complexity}, this framework allows for 
marginal constraints on the joint distributions. Furthermore, it incorporates the capacity for cost functions to depend on multiple marginals, a feature previously implied but now explicitly defined and integrated into our algorithm. 

\section{Multi-marginal optimal transport with global regulatization and probabilistic graphical models} \label{sec:global}

A line of recent works have drawn a connection between entropy-regularized multi-marginal optimal transport problems \eqref{eq:omt_multi_regularized} with a graph-structured cost \eqref{eq:cost_structure} and probabilistic graphical models (PGMs) \cite{HaaSinZha20,fan2021complexity,HaaRinChe20,SinHaaZha20,ringh2023mean,zhou2022efficient}.
In this section we review and generalize these connections and the computational methods for solving these types of optimal transport problems.

\subsection{A PGM formulation for multi-marginal optimal transport}

A PGM provides a compact representation for the joint distribution of a collection of random vectors, which have dependencies that can be described by a graph~\cite{KolFri09}. 
%
Define the tensors $K_\alpha(\bx_\alpha) = \exp (-C_\alpha(\bx_\alpha)/\epsilon)$ for $\alpha \in F$,
where $C_\alpha$ are the cost tensors in \cref{eq:cost_structure}, and $\epsilon$ is the regularization parameter in \cref{eq:omt_multi_regularized}.
Then, $\bK$ in \cref{eq:U} has the form $ \bK(\bx) = \prod_{\alpha\in F} K_\alpha(\bx_\alpha) $.
Moreover, to define the entropy regularization \cref{eq:entropy_term} for the graph-structured problem \cref{eq:ot_graph} we define the tensor $\bM$ as
$    \bM(\bx) = \prod_{\alpha\in F_\Gamma} \bmu_\alpha(\bx_\alpha)$.
The factorizations of the tensors $\bK$ and $\bM$ can be used to define the joint probability distribution of a PGM, 
\begin{equation} \label{eq:p_PGM}
    	p(\bx) = \frac{1}{Z}\left( \prod_{\alpha\in F} K_\alpha(\bx_\alpha) \right)   \left( \prod_{\alpha\in F_\Gamma} \bmu_\alpha(\bx_\alpha) \right)   ,
\end{equation}
where $Z$ is a normalization constant.
Moreover, note that the solution tensor of \eqref{eq:omt_multi_regularized} 
can be written as $\bB = \bK\odot\bU \odot \bM$ with
$\bU(\bx) = \prod_{\alpha\in F_\Gamma} \bu_\alpha(\bx_\alpha)$.
Thus, the scaling factors $\bu_\alpha$, for $\alpha\in F_\Gamma$ can be interpreted as local potentials of the modified graphical model $\bK(\bx)\bU(\bx) \bM(\bx)$ 
The Iterative Scaling algorithm \ref{alg:sinkhorn} finds potentials $\bu_\alpha$, for $\alpha\in F_\Gamma$, such that the graphical model defined by the tensor 
$\bB=\bK\odot\bU \odot \bM$
satisfies all the constraints $P_\alpha(\bB) = \bmu_\alpha$, for $\alpha\in F_\Gamma$.

\subsection{Iterative Scaling Belief Propagation}
Each step in the Iterative Scaling algorithm involves calculating a projection $P_j(\bK \odot \bU \odot \bM)$, which in terms of the PGM means inferring the distribution of the $j$-th variable.
This Bayesian inference problem can be addressed by the Belief Propagation algorithm~\cite[Part D]{mezard2009information}, \cite[Sec. 20.5.1]{KolFri09}.
Belief propagation is known to converge when the underlying graph is a tree, i.e., does not have any loops, and reduces the complexity of computing the projections in the Iterative Scaling algorithm immensely.
More precisely, whereas a brute force approach to compute the projections \cref{eq:proj_discrete} scales exponentially, Belief Propagation scales only linearly in the number of marginals \cite{HaaSinZha20}. \looseness=-1

In \cite{fan2021complexity}, a similar approach has be utilized for general graphical models with possible loops, by partitioning the graph as a so-called junction tree. A junction tree (also called tree decomposition) describes a partitioning of a graph, where several nodes are clustered
together, such that the interactions between the clusters can be described by a tree structure.
\begin{definition}[Junction tree \cite{dawid1992applications}] \label{def:junction_tree}
	A junction tree $\mathcal{T}= (\mathcal{C}, \mathcal{E})$ over a graph 
 {$G=(V,F)$}
 is a tree whose 
 nodes are cliques
 $c\in \mathcal{C}$ associated with subsets $\bx_c \subset V$, and 
 satisfying
 the 
 properties:
	\begin{itemize}
  \item Family preservation: For each $\alpha \in F$ there is a clique $ c \in \cC$ such that 
  $\alpha \subset c $.
\item Running intersection: For every pair of cliques $c_i, c_j \in \mathcal{C}$, every clique on the path between $c_i$ and $c_j$ contains $c_i \cap c_j$. 
	\end{itemize}
	For two adjoining cliques $c_i$ and $c_j$, define the separation set $  s_{ij}= \{ v\in V : v\in c_i \cap c_j \} $. 
\end{definition}

 A graph typically has several possible junction trees. In particular, one trivial junction tree is always given by a single cluster containing all nodes.
For our computational method, we set up the junction tree such that each constraint marginal or clique is represented by a leaf in the junction tree. 
The entropy-regularized graph structured optimal transport problem with junction tree decomposition $\mathcal{T}= (\mathcal{C}, \mathcal{E})$ is thus of the form 
\begin{equation} \label{eq:ot_jt}
		\min_{\bB \in \mR_+^{ d\times \dots \times d
  }} \langle \bC, \bB \rangle      + \epsilon \cH(\bB 
  {| \bM}  ) 
	\qquad	\text{ subject to }                               P_c (\bB) = \boldsymbol\mu_c,  \text { for } c \in \mathcal{C}_\Gamma,
\end{equation}
where $\mathcal{C}_\Gamma \subset \mathcal{C}$ is the set of constrained cliques,
and the cost tensor decouples as
\begin{equation}\label{eq:cost_jt}
	\bC(\bx) = \sum_{c \in \mathcal{C}} \bC_{c}(\bx_c).
\end{equation}
Such a junction tree representation exists for any graph, however it may require adding artificial clusters to the junction tree decomposition.
\begin{example}
    Some examples of valid junction trees for our method are shown in \cref{fig:traffic_graph_jt} (for the Euler flow example in \cref{fig:euler_graph}), \cref{fig:W2leastsquare_jt} (for the Wasserstein least squares problem \cref{eq:W2leastsquare_multimarginal} in \cref{fig:wasserstein_least_square}) , and \cref{fig:spline_jt} (for the Wasserstein splines problem in \cref{fig:wasserstein_spline}). 
\end{example}

\begin{figure}[tb]
	\caption{Junction tree for the graph from the Euler flow problem in 
 \cref{fig:euler_graph}, with leaf nodes corresponding to the constrained marginals.
 }
\centering
\scalebox{0.5}{
\begin{tikzpicture}[node distance=0.8cm] 
    \tikzstyle{circ}=[circle, minimum size = 10mm, thick, draw =black!80, node distance = 7mm] 
    \tikzstyle{circ_grey}=[circle, minimum size = 10mm, thick, draw =black!80, fill=gray!50, node distance = 7mm] 
    \tikzstyle{rect}=[rectangle, minimum size = 5mm, thick, draw =black!80, node distance = 7mm] 

    \node[circ] (c1) {$1,2,3$};
    \node[rect] (12) [right=of c1] {$1,3$};
    \node[circ] (c2) [right=of 12] {$1,3,4$};
    \node[rect] (23) [right=of c2] {$1,4$};
    \node[circ] (c3) [right=of 23] {$1,4,5$};
    \node[] (34) [right=of c3] {};
    \node[circ] (c4) [right=of 34] {$1,J-1,J$};
    \node[rect] (45) [right=of c4] {$1,J$};
    \node[circ_grey] (c5) [right=of 45] {$1,J$};
    
    \draw (c1) -- (12);
    \draw (12) -- (c2);
    \draw (c2) -- (23);
    \draw (23) -- (c3);
    \draw[dotted] (c3) -- (c4);
    \draw (c4) -- (45);
    \draw (45) -- (c5);
    
    \node[rect] (f2) [below=5mm of c1] {$2$}; 
    \node[circ_grey] (2) [below=5mm of f2] {$2$};
    \node[rect] (f3) [below=5mm of c2] {$3$}; 
    \node[circ_grey] (3) [below=5mm of f3] {$3$};
    \node[rect] (f4) [below=5mm of c3] {$4$}; 
    \node[circ_grey] (4) [below=5mm of f4] {$4$};
    \node[rect] (f5) [below=3mm of c4] {$J-1$}; 
    \node[circ_grey] (5) [below=3mm of f5] {$J-1$};
    
    \draw (c1) -- (f2); 
    \draw (f2) -- (2);
    \draw (c2) -- (f3); 
    \draw (f3) -- (3);
    \draw (c3) -- (f4); 
    \draw (f4) -- (4);
    \draw (c4) -- (f5); 
    \draw (f5) -- (5);
\end{tikzpicture} 
}
	\label{fig:traffic_graph_jt}
\end{figure}

\begin{figure}[tb]
	\centering
 	\caption{  Junction tree for the graph for the Wasserstein least square problem from \cref{fig:wasserstein_least_square}. The junction tree is constructed such that the leaves of the tree correspond to the constraints in the least square problem \cref{eq:W2leastsquare_multimarginal}. 
 }
	\label{fig:W2leastsquare_jt}
\scalebox{0.5}{
	\begin{tikzpicture}
		\tikzstyle{circ}=[circle, minimum size = 18mm, thick, draw =black!80, node distance = 8mm]
		\tikzstyle{rect}=[rectangle, minimum size = 10mm, thick, draw =black!80, node distance = 6mm]
		\tikzstyle{circ_s}=[circle, minimum size = 6mm, thick, draw =black!80, fill=gray!50, node distance = 6mm]
		\tikzstyle{rect_s}=[rectangle, minimum size = 6mm, thick, draw =black!80, node distance = 6mm]
		\node[circ] (c1) { \Large $\substack{0,\ 1,\\ J+1}$};
		\node[rect] (12) [right=of c1]{$0,\ J\!+\!1$};
		\node[circ] (c2) [right=of 12] {\Large $\substack{0,\ 2,\\ J+1}$};
		\node[rect] (23) [right=of c2] {$0,\ J\!+\!1$};
		\node[] (c3) [right=of 23] {};
		\node[circ] (cN) [right=of c3] {\Large $\substack{0,\ J-1,\\ J+1}$};
		\node[rect] (NNp1) [right=of cN] {$0,\ J\!+\!1$};
		\node[circ] (cNp1) [right=of NNp1] {\Large $\substack{0,\ J \\ J+1}$};
		\draw (c1) -- (12);   \draw (12) -- (c2);   \draw (c2) -- (23);
		\draw[dotted] (23) -- (cN);
		\draw (cN) -- (NNp1);	\draw (NNp1) -- (cNp1);
		\node[rect_s] (f1) [below=of c1] {$1$}; \node[circ_s] (1) [below =of f1] {$1$};
		\node[rect_s] (f2) [below=of c2] {$2$};	\node[circ_s] (2) [below=of f2] {$2$};
		\node[rect_s] (fNm1) [below=of cN] {$\!J\!\!-\!\!1\!$};	\node[circ_s] (Nm1) [below=of fNm1] {$\!\! J \!\!-\!\!1\!\!$};
		\node[rect_s] (fN) [below=of cNp1] {$J$}; \node[circ_s] (N) [below=of fN] {$J$};
		\draw (c1) -- (f1); \draw (f1) -- (1);
		\draw (c2) -- (f2); \draw (f2) -- (2);
		\draw (cN) -- (fNm1); \draw (fNm1) -- (Nm1);
		\draw (cNp1) -- (fN); \draw (fN) -- (N);
	\end{tikzpicture}
}
\end{figure}

\begin{figure}[tb]
	\centering
\caption{Junction tree of Wasserstein splines in \cref{fig:wasserstein_spline}. 
 }
	\label{fig:spline_jt}
 \scalebox{0.75}{
	\begin{tikzpicture}
		\tikzstyle{circ}=[circle, minimum size = 16mm, thick, draw =black!80, node distance = 3mm]
		\tikzstyle{scirc}=[circle, minimum size = 8mm, thick, draw =black!80, node distance = 3mm ,fill=gray!50,]  
		\tikzstyle{rect}=[rectangle, minimum size = 8mm, thick, draw =black!80, node distance = 4mm]
		\node[circ] (c1) {$\substack{v_0,v_1 \\ x_0,x_1}$};
		\node[rect] (s1) [right=of c1] {$\substack{v_1 \\ x_1}$};
		\node[circ] (c2) [right=of s1] {$\substack{v_1,v_2 \\ x_1,x_2}$};
		\node (s2) [right=of c2] {};
		\node[circ] (c3) [right=of s2] {$\substack{v_{J-2} \\ v_{J-1} \\ x_{J-2} \\ x_{J-1}}$};
		\node[rect] (s3) [right=of c3] {$\substack{v_{J-1} \\ x_{J-1}}$};
		\node[circ] (c4) [right=of s3] {$\substack{v_{J-1},v_{J} \\ x_{J-1},x_{J}}$};

		\node[rect] (s4) [below=of c1] {$x_0$};
		\node[rect] (s5) [below=of c2] {$x_{1}$};
		\node[rect] (s6) [below=of c3] {$x_{J\!-\!2\!}$};
		\node[rect] (s7) [below=of c4] {$x_{J\!-\!1\!}, x_{J}$};

		\draw (c1) -- (s1);
		\draw (c2) -- (s1);
		\draw[dotted] (c2) -- (c3);
		\draw (c3) -- (s3);
		\draw (c4) -- (s3);
		\draw (c1) -- (s4);
		\draw (s5) -- (c2);
		\draw (s6) -- (c3);
		\draw (s7) -- (c4);

		\node[scirc] (f1) [below=4mm of s4] {$x_0$};
        \node[scirc] (f2) [below=4mm of s5] {$x_1$};
        \node[scirc] (f3) [below=of s6] {\footnotesize $x_{J\!-\!2}$};
        \node[scirc] (f4) [below=of s7] {$\substack{x_{J -1} \\ x_J}$};

		\draw (s4) -- (f1);
        \draw (s5) -- (f2);
        \draw (s6) -- (f3);
        \draw (s7) -- (f4);
	\end{tikzpicture}
    }
\end{figure}

\begin{algorithm}[tb]
	\caption{Iterative Scaling Belief Propagation Algorithm for graphical 
 MOT
 with global regularization }
	\label{alg:isbp_jt}
	\KwIn{	Given a junction tree $\mathcal{T}= (\mathcal{C}, \mathcal{E})$ with constraints on the leaf cliques $\mathcal{C}_\Gamma \subset \mathcal{C}$}
		Initialize messages $m_{i\to j} (s_{ij}), \ m_{j\to i} (s_{ji}), \ \forall (i,j)\in \mathcal{E}$
   
   \While{not converged}{
   \vspace{-0.3cm}
     \begin{align}
		m_{i\to j} (s_{ij}) & \leftarrow \sum_{x_{c_i} \setminus s_{ij}} \bK_{c_i}(\bx_{c_i}) 
  \prod_{\ell \in N(i)\setminus j} m_{\ell\to i} (s_{\ell i}), \qquad  c_i \in \cC \setminus \cC_\Gamma \\
		m_{i \to j }(s_{ij}) & \leftarrow \mu_{c_i}(\bx_{c_i}) /  m_{j \to i} (\bx_{c_i}),  \qquad c_i \in \cC_\Gamma. 
\end{align}
}
\Return{$m_{i\to j} (s_{ij}),\ m_{j\to i} (s_{ji}), \ \forall (i,j)\in \mathcal{E}$}
\end{algorithm}

Similarly to the factorization \cref{eq:p_PGM}, we can then write the probability distribution over the PGM following the junction tree decomposition $\mathcal{T}= (\mathcal{C}, \mathcal{E})$ as a product
$
p( \bx) = \frac{1}{Z} \left( \prod_{c\in \mathcal{C}} \bK_c(\bx_c) \right) 
 \left( \prod_{c\in \mathcal{C}_\Gamma} \bmu_c(\bx_c) \right). 
 $
Let $N(i)$ be the neighboring nodes of the node $i$.
By extending the Iterative Scaling Belief Propagation methods in \cite{HaaSinZha20, fan2021complexity}, we propose \cref{alg:isbp_jt} to solve \cref{eq:ot_jt}.
Upon convergence of \cref{alg:isbp_jt}, the optimal transport plans on the clusters are given by 
\begin{equation} 
	\bB_c(\bx_c) =
 \bK_c(\bx_c) 
 \prod_{\ell \in N(i)} m_{\ell \to i} (s_{\ell i}).
\end{equation}
The algorithm sends messages between neighboring cliques. The messages are thus of the same size as the separation sets.
Recall that the leaves of the junction tree are associated with constraints in the optimal transport problem. Thus, the messages sent from these leaves correspond to an update step similar to the one in \cref{alg:sinkhorn}.
The messages sent within the graph encode the computation of the projections. 
More precisely, picking one leaf $c$ and sending the messages from all other leaves across the separation sets towards our leaf amounts to one exact computation of the projection 
$P_c(\bK \odot \bU \odot \bM)$.
However, the algorithm converges with any message schedule \cite{HaaSinZha20, fan2021complexity}.

Note that in the case where all cost are pair-wise and constraints act on one-marginals of the transport tensors, \cref{alg:isbp_jt} boils down to the Iterative Scaling Belief Propagation algorithms in \cite{HaaSinZha20, fan2021complexity}.
Here, we have extend on these methods and provide a more general algorithm for the case where we may have constraints on joint distributions, and where cost terms can depend on more than two marginals.


The Iterative Scaling Belief Propagation algorithm utilizes the graphical structure in problem \cref{eq:ot_jt}-\cref{eq:cost_jt} to significantly improve on the complexity of the Sinkhorn iterations.
Recall that the computational bottleneck of \cref{alg:sinkhorn} lies in the computation of the projections, which require $\mathcal{O}(d^J)$ operations.
\cref{alg:isbp_jt} reduces this complexity by utilizing the graphical structure for computing the projections. 
The complexity of the problem depends then on the size of the messages, i.e., the size of the separation sets. This quantity can also be described by the so-called tree-width.
\begin{definition}[Tree width \cite{halin1976s}]
	For a junction tree $\mathcal{T}= (\mathcal{C}, \mathcal{E})$, we define its width as $\text{width}(\cT) = \max_{c\in \cC} |c| -1$.
	For a graph $G$, we define its tree-width as $w(G) = \min\{ \text{width}(\cT) \ |\ \cT \text{ is a junction tree for } G\}$.
\end{definition}

In order for \cref{alg:isbp_jt} to be efficient we thus would like to pick a junction tree factorization with small tree-width. Given a graph, a junction tree decomposition that minimizes the tree-width is also called a minimal junction tree. 
In general, finding a junction tree with minimal width is a NP hard problem \cite{KolFri09}.
However, for all applications in this paper it is easy to find minimal junction trees.
\begin{example}\label{ex:wl}
The junction trees in Figures~\ref{fig:traffic_graph_jt}, \ref{fig:W2leastsquare_jt}, and \ref{fig:spline_jt} are minimal junction trees for the respective graphs.
Since the cliques in the former two junction trees are of cardinality 3, the respective graphs' tree-widths are 2.
Similarly, the last graph's tree width is 3.
\end{example}

For \cref{alg:isbp_jt}, we can directly apply the complexity analysis in \cite{fan2021complexity}.
We will discuss the complexity result in \cref{rem:complexity}.

\end{document}

%% file: ex_shared.tex
\usepackage{hyperref}
\usepackage{xr-hyper}

\usepackage{lipsum}
\usepackage{amsfonts}
\usepackage{graphicx}
\usepackage{epstopdf}
\ifpdf
  \DeclareGraphicsExtensions{.eps,.pdf,.png,.jpg}
\else
  \DeclareGraphicsExtensions{.eps}
\fi

\newsiamremark{remark}{Remark}
\newsiamremark{hypothesis}{Hypothesis}
\crefname{hypothesis}{Hypothesis}{Hypotheses}
\newsiamthm{claim}{Claim}
\newsiamremark{fact}{Fact}
\crefname{fact}{Fact}{Facts}

\headers{A parallel framework for graphical optimal transport}{J. Fan, I. Haasler, Q. Zhang, J. Karlsson, Y. Chen}

\title{A parallel framework for graphical optimal transport\thanks{ 
\funding{This work was supported by the Knut and Alice Wallenberg
foundation under grant KAW 2021.0274, the Swedish Research Council (VR) under grant 2020-03454, and the NSF under grant 2206576.}}}

\author{Jiaojiao Fan\thanks{NVIDIA,
\email{jiaojiaof@nvidia.com}. The first two authors contributed equally.}
       \and
 Isabel Haasler\thanks{Department of Information Technology, Uppsala University, Uppsala, Sweden, \email{isabel.haasler@it.uu.se }}
       \and
Qinsheng Zhang\thanks{NVIDIA, \email{qinshengz@nvidia.com} }
       \and
Johan Karlsson\thanks{Department of Mathematics,
       KTH Royal Institute of Technology,
       Stockholm, Sweden,
 \email{johan.karlsson@math.kth.se} }
       \and
Yongxin Chen\thanks{School of Aerospace Engineering,
       Georgia Institute of Technology,
       Atlanta, Georgia, USA   \email{yongchen@gatech.edu} }
       }

\usepackage{amsopn}
\DeclareMathOperator{\diag}{diag}

\usepackage{bm}
\usepackage{tikz}
\usetikzlibrary{positioning,chains,fit,shapes,calc}
\usepackage{algpseudocode}

\usepackage{subcaption}
\usepackage{enumerate}

\usepackage{verbatim}           %
\usepackage{xspace}             %
\usepackage{graphicx,float}     %
\usepackage{ifthen,calc}        %
\usepackage{textcomp}           %
\usepackage{fancybox}           %
\usepackage{hhline}             %
\usepackage{float}              %

\usepackage{booktabs}       
\usepackage{amsfonts}       
\usepackage{nicefrac}       
\usepackage{microtype}      

\usepackage{xr-hyper}

\usepackage{multirow}
\usepackage{bbm}
\usepackage{epstopdf}
\usepackage{algorithm,algcompatible}
\usepackage{soul}
\usepackage{tikz}
\usepackage{setspace}
\usepackage{enumitem}
\setenumerate[1]{label={(\arabic*)}}
\usepackage[algo2e, ruled,vlined, boxed]{algorithm2e}

\newcommand{\ba}{{\mathbf a}}
\newcommand{\bk}{{\mathbf k}}

\newcommand{\bu}{{\mathbf u}}
\newcommand{\bv}{{\mathbf v}}

\newcommand{\bg}{{\boldsymbol{\gamma}}}

\newcommand{\bq}{{\mathbf q}}

\newcommand{\bB}{{\mathbf B}}

\newcommand{\bx}{{\bf x}}
\def\bC{{\bf C}}
\def\bK{{\bf K}}
\def\bM{{\bf M}}
\def\bU{{\bf U}}
\newcommand{\bmu}{\boldsymbol{\mu}}
\newcommand{\ett}{{\bf 1}}
\def\minwrt[#1]{\underset{#1}{\text{minimize}}}
\def\maxwrt[#1]{\underset{#1}{\text{maximize }}}

\newcommand{\mR}{{\mathbb R}}

\newcommand{\cC}{{\mathcal C}}

\newcommand{\cH}{{\mathcal H}}

\newcommand{\cO}{{\mathcal O}}

\newcommand{\cQ}{{\mathcal Q}}

\newcommand{\cT}{{\mathcal T}}

\newcommand{\tr}{\operatorname{trace}}

\def\<{{\langle}}
\def\>{{\rangle}}
\def\holder{{H$\ddot{\text{o}}$lder's }}

\def\d{{\rm{d}}}
